\documentclass[11pt, a4paper]{article}

\usepackage[utf8]{inputenc}

\usepackage{color}

\usepackage{amsthm}

\usepackage{amsmath}
\usepackage{amssymb}
\usepackage{amstext}
\usepackage{amsfonts}
\usepackage{enumerate}

 \relax

\title{The algebraic semantics for the one-variable monadic fragment of the predicate logic ${\cal G}\forall_\sim$}

\author{\sc \scriptsize{Diego Casta\~no, Valeria Casta\~no, José Patricio Díaz Varela and Marcela Mu\~noz Santis} }

\date{}

\newtheorem{Theorem}{Theorem}[section] 
\newtheorem{Proposition}[Theorem]{Proposition} 
\newtheorem{Lemma}[Theorem]{Lemma} 
 
\newtheorem{Corollary}[Theorem]{Corollary} 
\newenvironment{Proof}{\noindent \bf Proof. \rm}{$\quad \hfill \blacksquare$ \medskip}  
\theoremstyle{definition}
\newtheorem{Definition}[Theorem]{ Definition} 
\newtheorem{Example}[Theorem]{Example} 
\newtheorem{Remark}[Theorem]{Remark} 

\newcommand{\s}{\mathop{\sim}}

\usepackage{pictex}

\begin{document}

\maketitle

\vspace{-1.8cm}

\begin{flushright}

\

\end{flushright}


\medskip

\begin{abstract}

In this article we characterize the equivalent algebraic semantics for  the one-variable monadic fragment of the first-order logic ${\cal G} \forall_{\s}$ defined by F. Esteva, L. Godo, P. Hájek and M. Navara in \textit{Residuated fuzzy logics with an involutive negation}, Archive for Mathematical Logic {\bf 39} (2000). To this end, we first introduce the variety  $\mathbb{MG}_{\s}$ as a certain class of Gödel algebras endowed with two monadic operators  and  a De Morgan negation. We study its basic properties, determine its subdirectly irreducible members and prove that this variety has the finite embeddabilty property. In particular, we prove that a special subvariety $\mathbb{CMG}_{\s}$ of $\mathbb{MG}_\sim$ is exactly the desired equivalent algebraic semantics; this is done via a functional representation of finite subdirectly irreducible algebras.

\end{abstract}

\section{Introduction and preliminares}

 $BL$-algebras were introduced by P. Hájek in his book \cite{Hj} as the algebraic semantics of his basic fuzzy logic (${\cal BL}$). One of the most important extensions of ${\cal BL}$ is Gödel Logic (${\cal G}$), obtained by adding the axiom schema $\varphi \to \varphi^2$, whose algebraic semantics, Gödel algebras,  is also given in \cite{Hj}. Recall that the variety of Gödel  algebras ($\mathbb{G}$) coincides with the subvariety of Heyting algebras determined by the prelineality equation $(x \to y) \vee (y \to x) \approx 1$ and $\mathbb{G}$ is generated by the Gödel t-norm $[0,1]_\mathbb{G}$.

Later, in \cite{EG}, F. Esteva, L, Godo, P. Hájek and M. Navara introduced the logic  ${\cal G}_{\s}$ as an expansion of Gödel logic adding an arbitrary independent involutive negation given by  the following axiom schemata:

\

\begin{minipage}{8cm}
\begin{enumerate}[$(\s 1)$]

\item $\s(\s \varphi) \equiv \varphi$

\item $\neg \varphi \to \s \varphi$

\item  $\Delta (\varphi \to \psi) \to \Delta (\s \psi \to \s \varphi)$

\end{enumerate}
\end{minipage}
\begin{minipage}{6cm}
\begin{enumerate}[$(\Delta 1)$]
\item $\Delta \varphi \vee \neg \Delta \varphi$

\item $\Delta(\varphi \vee \psi) \to (\Delta \varphi \vee \Delta \psi)$

\item[$(\Delta 5)$] $\Delta (\varphi \to \psi)\to (\Delta \varphi \to \Delta \psi)$

\end{enumerate}
\end{minipage}

\

\noindent where $\Delta \varphi$ is $\neg \s \varphi$ and $\varphi \equiv \psi$ is $(\varphi \to \psi) \& (\psi \to \varphi)$. The deduction rules of ${\cal G}_{\s}$ are {\it modus ponens}  ($\displaystyle \frac{\varphi, \varphi \to \psi}{\psi}$) and {\it necessitation} for $\Delta$ ($\displaystyle \frac{\varphi}{\Delta \varphi}$).

  They  also defined the  predicate calculus based on  this logic $({\cal G} \forall_{\s})$ by taking as axioms those of ${\cal G}_{\s}$ plus the following five axioms for quantifiers given by  P. Hájek in his basic predicate logic:

\begin{enumerate}

\item[$(\forall 1)$] $(\forall x) \varphi(x) \to \varphi(t)$ \ ($t$ substitutable for $x$ in $\varphi(x)$) 

\item[$(\exists 1)$] $ \varphi(t) \to (\exists x)\varphi(x)$ \ ($t$ substitutable for $x$ in $\varphi(x)$)

\item[$(\forall 2)$] $(\forall x) (\psi \to \varphi) \to (\psi \to (\forall x) \varphi) $ \ ($x$ not free in $\psi$) 

\item[$(\exists 2)$] $(\forall x) (\psi \to \varphi) \to ((\exists x) \psi \to  \varphi) $ \ ($x$ not free in $\varphi$)

\item[$(\forall 3)$] $(\forall x) (\psi \vee \varphi) \to ((\forall x) \psi \vee \varphi) $ \ ($x$ not free in $\varphi$)

\end{enumerate}

\noindent The inference rules of  this calculus are modus ponens, necessitation ($\displaystyle \frac{\varphi}{\Delta \varphi}$) and generalization ($\displaystyle \frac{\varphi}{\forall \varphi}$).

In  recent articles, D. Castaño, C. Cimadamore, J. P. Díaz Varela and L. Rueda introduced and studied the variety $\mathbb{MG}$ of monadic Gödel algebras as an expansion of Gödel algebras by two monadic operators $\exists$ and $\forall$. They showed that these algebras constitute the equivalent algebraic semantics of the one-variable monadic fragment of the predicate  Gödel logic \cite{Bahia, Complet, Laura}. Following this line of research we intend to find the equivalent algebraic semantics of the one-variable monadic fragment of the first-order logic ${\cal G}\forall_{\s}$, which is the main result of this paper. 
We start by introducing  an expansion $\mathbb{MG}_{\s}$ of the variety $\mathbb{MG}$  by a unary operation  $\s$. In Sections 2 and 3, we make a standard study of the basic properties of the algebras in this variety, which includes the  characterization of  congruences and subdirectly irreducible algebras. We also prove that $\mathbb{MG}_\sim$ has the finite embeddability property (FEP for short) and so it  is generated by its finite members. In the next sections,  we  introduce an important subvariety of $\mathbb{MG}_{\s}$, the variety $\mathbb{CMG}_{\s}$, characterized within the variety $\mathbb{MG}_{\s}$ by the identity $\exists (x \wedge \s x) \to  \forall (x \vee \s x) \approx 1 $, that is, $\exists (x \wedge \s x) \leq  \forall (x \vee \s x)$. We are interested in this subvariety because  the formula $(\exists x) (\varphi \wedge \s \varphi) \to (\forall x) (\varphi \vee \s \varphi)$ is a theorem in the predicate calculus based on the logic ${\cal G}_{\s}$. We give a complete characterization of the range of the monadic operators in algebras within this variety, which  allows us to show that $\mathbb{CMG}_{\s}$ is precisely the variety generated by   the functional monadic $G_{\s}$-algebras. Finally, in Section 6, taking into account that functional algebras encode the information of the models used to interpret the one-variable monadic fragment of the first-order logic ${\cal G}\forall_{\s}$, we prove that  $\mathbb{CMG}_{\s}$ is the equivalent algebraic semantics for this calculus.

\

Recall (see \cite{Bahia}) that the algebras in $\mathbb{MG}$ (monadic Gödel algebras, which we call $MG$-algebras in this article, for short) are Gödel algebras endowed with two unary operations $\forall$ and $\exists$ that satisfy the following identities:

  \begin{enumerate}[(M1)]
\begin{minipage}{7cm}
   \item $\forall x\to x\thickapprox 1,$
    \item $\forall(x\to\forall y)\thickapprox\exists x\to\forall y,$
    
\end{minipage}
\begin{minipage}{7cm}    
    
\item $\forall(\forall x\to y) \thickapprox \forall x\to\forall y,$    
    \item $\forall(\exists x\lor y) \thickapprox \exists x\lor\forall y.$
\end{minipage}
  \end{enumerate}
Note that in \cite{Bahia} there is an additional identity (M5)  that holds trivially  in this case.

\

We are interested in studying  $MG$-algebras with an involutive negation. The following definition is motivated by the axioms of the predicate calculus based on the logic ${\cal G}_{\s}$.

\begin{Definition}
An algebra  $\langle A, \vee, \wedge, \to, \s, \exists, \forall, 0, 1 \rangle$ is an {\it $MG_{\s}$-algebra} if $\langle A, \vee, \wedge, \to, \exists, \forall, 0, 1 \rangle $ is a monadic Gödel algebra, $\langle A, \vee, \wedge, \s, 0, 1 \rangle$ is a  De Morgan algebra and the following identities hold:
\begin{enumerate}

\item[$ (N)$] $ \neg x\to \s x \approx 1$,  that is,  $ \neg x\leq \s x$,

\item[  $(Q)$] $  \forall x \approx \s \exists\s x$,

\end{enumerate}
 where $\neg x := x \to 0$.
 
\end{Definition}

We denote by $\mathbb{MG}_{\s}$ the variety of  $ MG_{\s}$-algebras and by $\mathbb{G}_{\s}$ the variety of Gödel algebras with a De Morgan negation that satisfies $(N)$.

Here and subsequently ${\bf A}$ denotes a Gödel algebra, and $\langle {\bf A}, \s, \exists, \forall \rangle$ the Gödel algebra ${\bf A}$ enriched with a monadic structure and with an involutive negation. Observe that if $\langle {\bf A}, \s, \exists, \forall \rangle$ is an $MG_{\s}$-algebra, $\langle {\bf A}, \exists, \forall \rangle$ is an $MG$-algebra and $\langle {\bf A},  \s \rangle$ is a $G_{\s}$-algebra. For this reason we are interested in the properties of the algebras in $\mathbb{MG}$ and $\mathbb{G}_{\s}$. To close this section we summarize some of the main properties of $G_{\s}$-algebras and $MG$-algebras which will be useful in the development of  this article.

The variety $\mathbb{G}_{\s}$ was introduced by A. Monteiro  in \cite{Mont} as the class of symmetric fully linear algebras and later studied in, for example,  \cite{Pato}. The proofs of the following properties can be found in \cite{Mont} and \cite{Pato}. We start by recalling some  arithmetical properties.

\begin{Lemma} \label{prop aritmeticas de $G_sim$}
  Let $\langle {\bf A}, \s \rangle$ be an algebra in $\mathbb{G}_{\s}$ and $a,b\in A$. Then:

\begin{enumerate}[$(1)$]
\begin{minipage}{8cm}
 
    \item $\neg a \vee \neg \neg a =1$,
    \item $\Delta (a\to b)=\Delta (\s b\to \s a)$,
    \item $\Delta a\lor \neg\Delta a=1$,        \end{minipage}
    \begin{minipage}{5cm}
    \item $\Delta(a\lor b)= \Delta a\lor \Delta b$, 
    \item $\Delta a\land \Delta (a\to b)\leq \Delta b$,
    \item $\Delta a \leq a$; \  $\Delta \Delta a = \Delta a$,
\end{minipage}
\end{enumerate}

\noindent where $\Delta a:=\neg\s a$.

\end{Lemma}

Observe that $(1)$ says that algebras in $\mathbb{G}_{\s}$ are  Stone algebras and  (2)-(5) say that they are $SBL_{\s}$-algebras (see \cite{EG}). 

Recall that a {\it filter} of a $G_{\s}$-algebra $\langle {\bf A}, \s \rangle$ is  an increasing nonempty set $F \subseteq A$ such that $a \wedge b \in F$ whenever $a,b \in F$. A filter $F$ is called {\it prime} if $a  \vee b \in F$ implies $a \in F$ or $b \in F$. We have the following important properties of the variety $\mathbb{G}_{\s}$.

\begin{Lemma} \label{propiedades de Godel simetricas}

Let $\langle {\bf A}, \s \rangle$ be an algebra in $\mathbb{G}_{\s}$. Then:

\begin{enumerate}[$(1)$]

\item The condition $(N)$ 
in the definition of $\mathbb{G}_{\s}$ can be replaced by the Kleene condition $(K)$: $x \wedge \s x \leq y \vee \s y$.

\item The set of prime filters of $\langle {\bf A}, \s \rangle$ is a disjoint union of chains. Hence, if $\langle {\bf A}, \s \rangle$ is a finite algebra, $\langle {\bf A}, \s \rangle$ is a product of totally ordered algebras.

\item The congruence lattice  and the lattice of regular filters of $\langle {\bf A}, \s \rangle$ are isomorphic, where $F$ is a regular filter if $F$ is a filter that satisfies $\Delta x \in F$ whenever $x \in F$ (see \cite{Sanka}).

\item $\langle {\bf A}, \s \rangle$ is a subdirectly irreducible algebra if and only if $\langle {\bf A}, \s \rangle$ is totally ordered if and only if $\langle {\bf A}, \s \rangle$ is simple.

\end{enumerate}

\end{Lemma}

Next we summarize some of the basic properties of monadic Gödel algebras; all the proofs can be found in \cite{Bahia} in the broader context of monadic $BL$-algebras.

\begin{Lemma} \label{propiedades de las monádicas de Gödel}
Let $\langle {\bf A}, \exists, \forall \rangle$ be an algebra in $\mathbb{MG}$. Then:

\begin{enumerate}[$(1)$]
\item  $\exists A = \forall A$.

\item $\exists {\bf A}$ is a subalgebra of ${\bf A}$.

\item  The congruence lattice and the lattice of monadic filters of $\langle {\bf A}, \exists, \forall \rangle$ are isomorphic, where $F$ is a monadic filter if $F$ is a filter that satisfies $\forall x \in F$ whenever $x \in F$.

\item The congruence lattices of $\langle {\bf A}, \exists, \forall \rangle$ and $\exists {\bf A}$ are isomorphic.

\item $\langle {\bf A}, \exists, \forall \rangle$ is finitely subdirectly irreducible if and only if $\exists {\bf A}$ is totally ordered.

\end{enumerate}
\end{Lemma}

Some arithmetical properties of $MG$-algebras useful for this work are the following.

\begin{Lemma} \label{Propiedades algebraicas de las monádicas}
Let $\langle {\bf A}, \exists, \forall \rangle$ be an algebra in $\mathbb{MG}$. Then, for any $a,b \in A$:

\begin{enumerate}[$(1)$]
\begin{minipage}{7cm}
\item $\forall 1 = \exists 1=1$ and $\forall 0 = \exists 0 =0$, 
\item $\forall a \leq a$; \ $a \leq \exists a$,
\item $\forall\forall a = \forall a$ and $\exists\exists a = \exists a$, 

\end{minipage}
\begin{minipage}{6cm}
\item $\forall (a \vee \forall b) = \forall a \vee \forall b$,
\item if $a \leq b$, then $\forall a \leq \forall b$ and $\exists a \leq \exists b$.
\end{minipage}
\end{enumerate}

\end{Lemma}

\section{Subdirectly irreducible $MG_{\s}$-algebras}

Let $\langle {\bf A}, \s, \exists, \forall \rangle $ be an algebra in the variety $\mathbb{MG}_{\s}$. In this section we prove that $\langle \exists {\bf A}, \s \rangle$ is a subalgebra of $\langle {\bf A}, \s \rangle$ and that the congruence lattices of $\langle {\bf A}, \s, \exists, \forall \rangle$ and $\langle \exists {\bf A}, \s \rangle$ are isomorphic.

\begin{Lemma} \label{el rango es subálebra de Godel simétrica}
If  $\langle {\bf A}, \s, \exists, \forall \rangle$ is  an algebra in $\mathbb{MG}_{\s}$, then  $\langle \exists {\bf A}, \s \rangle$ is a subalgebra of $\langle {\bf A}, \s \rangle$.  
\end{Lemma}

\begin{Proof}
Since $\langle {\bf A}, \exists, \forall \rangle$ is an $MG$-algebra,  by Lemma \ref{propiedades de las monádicas de Gödel}, we have that $\exists {\bf A}$ is a Gödel algebra and $\exists A= \forall A$. We now consider $x \in \exists A$, then $x=\exists z$, for some $z \in A$. From $(Q)$, we obtain $\forall \s z = \s \exists z$ and thus $\s x = \forall \s z$. This proves that $\langle \exists {\bf A}, \s \rangle $ is a subalgebra of $\langle {\bf A}, \s \rangle$.
\end{Proof}

Let us define the following lattices  

\begin{itemize}
\item $\mbox{{\bf Con}}({\bf A})$: the congruence lattice of an algebra ${\bf A}$.

\item  ${\bf F}_r(\langle {\bf A}, \s \rangle)$: the lattice of regular filters of an algebra $\langle {\bf A}, \s \rangle$ in $\mathbb{G}_{\s}$.

\item  ${\bf F}_m(\langle {\bf A}, \exists, \forall \rangle)$: the lattice of monadic filters of an algebra $\langle {\bf A}, \exists, \forall \rangle$ in $\mathbb{MG}$.  

\item   ${\bf F}_{mr}(\langle {\bf A}, \s, \exists, \forall \rangle)$: the lattice of monadic regular filters of an algebra $\langle {\bf A}, \s, \exists, \forall \rangle$ in $\mathbb{MG}_{\s}$, where $F \subseteq  A$ is a  {\it monadic regular filter} if $F$ is both, monadic and regular. 
\end{itemize}

In \cite[Th. 2.7]{Bahia} it was proved that ${\bf F}_m(\langle {\bf A}, \exists, \forall \rangle)$ is isomorphic to ${\bf Con}(\langle {\bf A}, \exists, \forall \rangle))$, for each algebra $\langle {\bf A}, \exists, \forall \rangle$ in $\mathbb{MG}$, and the isomorphism is given by $F \longmapsto \theta_F=\{(a,b) \in A^2: (a \to b) \wedge (b \to a) \in F\}.$

Similarly, in \cite[Th. 3.3]{Sanka} it was proved that ${\bf F}_r(\langle {\bf A}, \s \rangle)$ is isomorphic to ${\bf Con}(\langle {\bf A}, \s \rangle)$, for each algebra $\langle {\bf A}, \s \rangle$ in $\mathbb{G}_{\s}$, and the isomorphism is given by $F \longmapsto \theta'_F=\{(a,b) \in A^2: a \wedge f = b \wedge f, \mbox{ for some } f \in F\}.$

Taking these results into account as well as the fact that $\theta_F=\theta'_F$ in an $MG_{\s}$-algebra the following proposition is immediate.

\begin{Proposition}
Let $\langle {\bf A}, \s, \exists, \forall \rangle$ be an ${MG}_{\s}$-algebra. The correspondence  $   \mbox{{\bf F}}_{mr}(\langle {\bf A}, \s, \exists, \forall \rangle) \to \mbox{{\bf Con}}(\langle {\bf A}, \s, \exists, \forall \rangle)$ defined by  $F \longmapsto \theta_F=\{(a,b) \in A^2: (a \to b) \wedge (b \to a) \in F\}$  is an order isomorphism whose inverse is given by $\theta \longmapsto 1/\theta$.

\end{Proposition}

In $\cite{Sanka}$ H.P. Sankappanavar defined the following terms:

$$t_k(x):=\Delta^0 x \wedge \Delta^1 x \wedge \Delta^2 x \wedge \ldots \wedge \Delta^k x,$$

\noindent where $k \in \mathbb{N}_0$ and $\Delta^n x$ is defined recursively: $\Delta^0x:=x$ and $\Delta^{n+1}x:=\Delta^n(\Delta x)$, for $n \in \mathbb{N}_0$. He used them to  prove that given a De Morgan Heyting algebra ${\bf L}$ and $X \subseteq L$, $X \ne \emptyset$,  the regular filter generated by $X$ is the set 

$Fg_r(X) = \{ a \in L: t_{k_1}(x_1) \wedge \ldots \wedge t_{k_n}(x_n)\leq a, \mbox{ for some } x_1, \ldots,  x_n \in X,  k_i  \in \mathbb{N}_0, n \in \mathbb{N} \}.$

Observe that algebras in the variety $\mathbb{G}_{\s}$ satisfy $\Delta x \leq  x$ and $\Delta^2x=\Delta x$ (Lemma \ref{prop aritmeticas de $G_sim$}). Then, on any of these algebras we have $t_0(x)= x$ and  $t_k(x)= \Delta x $ for all $k \in \mathbb{N}$. Thus, for $\langle {\bf A}, \s \rangle$ in $\mathbb{G}_{\s}$ and $\emptyset \ne X \subseteq A$:

\begin{eqnarray} \label{filtro monádico regular}
Fg_r(X) = \{ a \in A: \Delta x_1 \wedge \ldots  \wedge \Delta x_n\leq a, \mbox{ for some } x_1, \ldots,  x_n \in X,   n \in \mathbb{N} \}.
\end{eqnarray}

\begin{Lemma} \label{filtro monadico regular generado}
Given an ${MG}_{\s}$-algebra $\langle {\bf A}, \s, \exists, \forall \rangle$ and $X \subseteq A$, $X \ne \emptyset$, the monadic regular filter generated by $X$ is the set

\noindent
$\begin{array}{lcl}
Fg_{mr}(X) & = & Fg_r(\{\forall x: x \in X\}) \\
& = & \{a \in A:  \Delta\forall x_1 \wedge  \cdots \wedge  \Delta\forall x_n  \leq a, \mbox{ where } x_1, x_2, \ldots, x_n \in X, n \in \mathbb{N}\}.

\end{array}$

\end{Lemma}
\begin{Proof}
Let $S=\{a \in A:  \Delta\forall x_1 \wedge  \cdots \wedge  \Delta\forall x_n  \leq a, \mbox{ where } x_1, x_2, \ldots, x_n \in X\}$. By (\ref{filtro monádico regular}) we have that $S=Fg_r(\{\forall x, x \in X\}$ in the algebra $\langle {\bf A}, \s \rangle$. Thus, if $a \in S$, $\Delta a \in S$. 

Let us prove that $\forall a \in S$ whenever $a \in S$. Given  $a \in S$, there exist $x_1, \ldots, x_n \in X$  such that  $ \Delta\forall x_1 \wedge  \cdots \wedge \Delta\forall x_n \leq a$.   By Lemma \ref{el rango es subálebra de Godel simétrica} and Lemma \ref{propiedades de las monádicas de Gödel} (1)  we know that $\langle \forall {\bf A}, \s \rangle$ is a subalgebra of $\langle {\bf A}, \s \rangle$.  Then, there exists $y \in A$ such that $\Delta\forall x_1 \wedge  \cdots \wedge  \Delta\forall x_n= \forall y$.  From this $\forall y \to a =1$ and so, $\forall(\forall y \to a)=\forall 1 =1$.  Finally, by $(M3)$ we obtain $\forall y \to \forall a = 1$, or equivalently  $ \Delta\forall x_1 \wedge  \cdots \wedge  \Delta\forall x_n \leq \forall a$.

We have proved that $S$ is a monadic regular filter, and it is clearly the  smallest one that contains $X$. 
\end{Proof}

\begin{Theorem} \label{congurncias de A y del rango coinciden}

Let $\langle {\bf A}, \s, \exists, \forall \rangle $ be an $MG_{\s}$-algebra. The correspondence ${\bf F}_{mr}(\langle {\bf A}, \s, \exists,  \forall \rangle) \to {\bf F}_r(\langle \exists {\bf A}, \s \rangle)$ defined by $F \longmapsto F \cap \exists A$ is an order isomorphism whose inverse is given by $S \mapsto Fg_{mr}(S)$.
\end{Theorem}

\begin{Proof}
Let $F$ be a monadic regular filter of $\langle {\bf A}, \s, \exists, \forall \rangle$. Since  $\langle \exists {\bf A}, \s \rangle$ is a subalgebra of $\langle {\bf A}, \s \rangle$, we have that $F \cap \exists A$ is a regular filter of $\langle \exists {\bf A}, \s \rangle$.

Let us prove that $Fg_{mr}(F \cap \exists A)=F$, for all $F \in \mbox{F}_{mr}(\langle {\bf A}, \s, \exists, \forall \rangle)$. From $F \cap \exists A \subseteq F$, it follows that $Fg_{mr}(F \cap \exists A)\subseteq F$. Now suppose $x \in F$. Since $\exists A=\forall A$ and $F$ is a monadic filter we have that $\forall x \in F \cap \exists A$. By Lemma \ref{Propiedades algebraicas de las monádicas}(2) we know that $ \forall x \leq x$ and consequently $x \in Fg_{mr}(F \cap \exists A)$.

Let us prove that $Fg_{mr}(S) \cap \exists A =S$ for each $S \in \mbox{F}_r(\langle \exists {\bf A}, \s \rangle)$. Using Lemma \ref{filtro monadico regular generado} and $\forall S = S$ we obtain that $$Fg_{mr}(S) \cap \exists A = Fg_{r}(\forall S) \cap \exists A = Fg_{r}(S) \cap \exists A.$$

Moreover, $Fg_{r}(S) \cap \exists A = S$. Indeed, from $(\ref{filtro monádico regular})$ we know that if $a \in Fg_r(S) \cap \exists A$ then  $a \geq \Delta x_1 \wedge \ldots \wedge \Delta x_n$, for some $x_1, \ldots, x_n \in S$. Now, since  $S$ is a regular filter, $a \in \exists A$ and $S$ is an increasing set on $\exists A$, we have that $a \in S$. Thus,  $Fg_{r}(S) \cap \exists A \subseteq S$.  The other inclusion is trivial. 

Hence, it is clear that the correspondence $F \mapsto F \cap \exists A$ is an order isomorphism.  
\end{Proof}

From Lemma  \ref{propiedades de Godel simetricas}   we obtain the following immediate results.

\begin{Corollary}\label{subdirectamente irreducibles MGalgebras simétricas} Let $\langle {\bf A}, \s, \exists, \forall \rangle \in \mathbb{MG}_{\s}$, then:

\begin{enumerate}[$(1)$]

 \item $\mbox{{\bf Con}}(\langle {\bf A}, \s, \exists, \forall \rangle) \cong {\bf F}_{mr}(\langle {\bf A}, \s, \exists, \forall \rangle) \cong {\bf F}_r(\langle \exists {\bf A}, \s \rangle) \cong {\bf Con}(\langle \exists {\bf A}, \s \rangle)$.

\item The following are equivalent:
\begin{enumerate}
\item[$(a)$] $\langle {\bf A}, \s, \exists, \forall \rangle$ is subdirectly irreducible. 
\item[$(b)$] $\langle {\bf A}, \s, \exists, \forall \rangle$ is simple.
\item[$(c)$]  $\langle \exists {\bf A}, \s \rangle$ is subdirectly irreducible.
\item[$(d)$]  $\langle \exists {\bf A}, \s \rangle$ is simple.
\item[$(e)$]  $\langle \exists {\bf A}, \s \rangle $ is totally ordered.
\end{enumerate}

\end{enumerate}

\end{Corollary}

\begin{Remark} \label{las si finitas son producto directo de cadenas}
In particular, if  $\langle {\bf A}, \s, \exists, \forall \rangle$ is a finite subdirectly irreducible algebra in $\mathbb{MG}_{\s}$ we have, by Corollary \ref{subdirectamente irreducibles MGalgebras simétricas} and  Lemma \ref{propiedades de Godel simetricas} (2), that $\langle {\bf A}, \s\rangle \cong \prod_{i=1}^r \langle {\bf C}_{n_i}, \s \rangle$, where $\langle {\bf C}_{n_i}, \s \rangle$ is the $G_{\s}$-chain with $n_i$ elements, and $\langle \exists {\bf A}, \s \rangle$ is a totally ordered subalgebra of $\langle {\bf A}, \s \rangle$.

\end{Remark}

As a consequence of the results given in this section we obtain an important property of the variety $\mathbb{MG}_{\s}$.
\begin{Theorem}

$\mathbb{MG}_{\s}$ is a discriminator variety. 
\end{Theorem}

\begin{Proof}
We consider the term $T(x):= \Delta \forall x$. Since $\langle \exists {\bf A}, \s \rangle$ is totally ordered for each subdirectly irreducible algebra  $\langle {\bf A}, \s, \exists, \forall \rangle$ in $\mathbb{MG}_{\s}$, it follows that:

$$T(a) = \left\{ \begin{array}{ll}
1 & \mbox{ if } a =1 \\
0 & \mbox{ if } a \ne 1,
\end{array} \right.$$
for all $a \in A$.
Then, it is immediate that a discriminator term for each subdirectly irreducible algebra is given by $$t(x,y,z) := (T((x \to y) \wedge (y \to x)) \wedge z) \vee (\neg (T((x \to y) \wedge (y \to x)) \wedge x).$$
\end{Proof}

\begin{Corollary} \label{si = simple=indescomponible}
 $\langle {\bf A}, \s \exists, \forall \rangle \in \mathbb{MG}_{\s}$ is subdirectly irreducible if and only if $\langle {\bf A}, \s, \exists, \forall \rangle $ is simple if and only if $\langle {\bf A}, \s,  \exists, \forall \rangle$ is directly indescomposable.
\end{Corollary}

\section{Generations by finite members}

We note that the variety $\mathbb{G}_{\s}$ is locally finite. In fact, since all of the subdirectly irreducible algebras in this variety are totally ordered (Lemma \ref{propiedades de Godel simetricas}) it is easy to see that each $n$-generated chain  has at most $2n+2$ elements. This implies that the class of subdirectly irreducible algebras in $\mathbb{G}_{\s}$ is uniformly locally finite and so, using \cite[Th. 3.7]{Bez} we have that $\mathbb{G}_{\s}$ is locally finite. This property is not true for the 
variety $\mathbb{MG}_{\s}$, as we will show later. 
Nevertheless,  we prove next that the variety $\mathbb{MG}_{\s}$ has the FEP and thus, it is generated by its finite members. In order to show this it will  be  useful to describe which subalgebras of a $G_\sim$-algebra can be the range  of the quantifiers.

First recall that, in \cite{Bahia}, D. Castaño et al.\ give a characterization of the range of the quantifiers of a monadic $BL$-algebra ${\bf L}$ by means of certain types of subalgebras  called {\it $m$-relatively complete subalgebras}. Given  a $BL$-algebra ${\bf L}$, a subalgebra ${\bf C}$ of ${\bf L}$ is {\it $m$-relatively complete} if it satisfies:

\begin{enumerate}[$(s_1)$]
\item For every  $a \in L$, the set $\{c \in C: c \leq a\}$ has a greatest element and $\{c \in C: c \geq a\}$ has a least element.

\item For every $a \in L$ and $c_1, c_2 \in C$ such that $c_1 \leq c_2 \vee a$, there exists $c_3 \in C$ such that $c_1 \leq c_2 \vee c_3$ and $c_3 \leq a$.

\item For every $a \in L$ and $c_1 \in C$ such that $a \ast a \leq c_1$, there exists $c_2 \in C$ such that $a \leq c_2$ and $c_2 \ast c_2 \leq c_1$.
\end{enumerate}

Given a $G_{\s}$-algebra $\langle {\bf A},  \s \rangle$,  a subalgebra $\langle {\bf C}, \s \rangle$  of $\langle {\bf A},  \s \rangle$    is {\it $m$-relatively complete}  if it satisfies $(s_1)$-$(s_3)$ from the above definition. In this case, since $\wedge$ coincides with $\ast$, the condition $(s_3)$ trivially holds. Moreover, if $\langle {\bf C}, \s \rangle$ is totally ordered, $(s_2)$ may  be replaced by the following simpler equivalent condition (see \cite{Bahia}):

\begin{enumerate}

\item[$(s_2')$] If $1=c \vee a$ for some $c \in C$ and $a \in A$, then $c=1$ or $a=1$.

\end{enumerate}

The following characterization of the range of the quantifier $\exists$ in an $MG_{\s}$-algebra is analogous to the result given in \cite{Bahia} for $MG$-algebras, and it is useful to build $MG_{\s}$-algebras from a given $G_{\s}$-algebra.

\begin{Theorem} \label{cuantificador asociado a C}
Let $\langle {\bf A}, \s \rangle$ be a ${G}_{\s}$-algebra  and $\langle {\bf C}, \s \rangle$ an  $m$-relatively complete subalgebra of $\langle {\bf A}, \s \rangle$. We define on $A$ the operations
$$\exists a:= \min\{c \in C: c\geq a\}, \ \ \ \ \forall a:= \max\{c \in C: c\leq a\},$$ then $\langle {\bf A}, \s, \exists, \forall \rangle$ is an  $MG_{\s}$-algebra such that $\langle \forall {\bf A}, \s \rangle = \langle \exists {\bf A}, \s \rangle =\langle {\bf C}, \s \rangle$. Conversely, if $\langle {\bf A}, \s, \exists, \forall \rangle$ is an $MG_{\s}$-algebra, then $\langle \forall {\bf A}, \s \rangle = \langle \exists {\bf A}, \s \rangle$ is an  $m$-relatively complete subalgebra of $\langle {\bf A}, \s \rangle$.
\end{Theorem}

\begin{Proof}
Suppose that $\langle {\bf C}, \s \rangle$ is an $m$-relatively complete subalgebra of $\langle {\bf A}, \s \rangle$ and $\forall$ and $\exists$ are defined as in the statement of the theorem. In \cite{Bahia} it was proved that $\langle {\bf A}, \exists, \forall \rangle$ is a monadic $BL$-algebra and consequently, it is an $MG$-algebra. Moreover, $\exists A = \forall A = C$. It only remains to show that $\forall a = \s \exists \s a $ for each $a \in A$.
We write $d =  \forall a = \max\{c \in C: c\leq a\}$ and $m= \exists \s a= \min\{c \in C: c\geq \s a\}$.
From $m \geq \s a$ we have that $\s m \leq a$. Since $\langle {\bf C}, \s \rangle$ is a subalgebra of $\langle {\bf A}, \s \rangle$ and $m \in C$ we obtain $\s m \in C$ and so, $\s m \leq d$. Dually, from  $d \leq  a$, $\s d \geq \s a$ and since $\s d \in C$, we have that $m \leq \s d$ and so $\s m \geq d$.
Therefore, $\langle {\bf A}, \s, \exists, \forall \rangle$ is an $MG_{\s}$-algebra such that $\langle \forall {\bf A}, \s \rangle = \langle \exists {\bf A}, \s \rangle =\langle {\bf C}, \s \rangle$.

The converse  follows from Lemma \ref{el rango es subálebra de Godel simétrica} and \cite[Th. 3.1]{Bahia}. 
\end{Proof}

\begin{Example} \label{ejemplos de construcción de MG algebras simetricas} 

Fix two natural numbers $n$ and  $r$, and consider the  chain $\langle {\bf C}_n, \s \rangle$ with $n$ elements in $\mathbb{G}_{\s}$ and the $G_{\s}$-algebra $\langle {\bf A}, \s \rangle:= \langle {\bf C}_n, \s \rangle^r$. It is easy to check that $\langle {\bf D}_1, \s \rangle$ and $\langle {\bf D}_2, \s \rangle$ are $m$-relatively  complete subalgebras of $\langle {\bf A}, \s \rangle$, where $D_1:=\{x \in A: x(i)=x(j), \mbox{ for all } i,j \in \{1, \ldots, n\}\}$, $D_2:= \{0, 1\}$ and $0$, $1$ are the least and greatest  elements of ${\bf A}$, respectively. Then, by  Theorem \ref{cuantificador asociado a C} we have that $\langle {\bf A}, \s,  \exists_1, \forall_1 \rangle $ and $\langle {\bf A}, \s, \exists_2, \forall_2 \rangle $ are both $MG_{\s}$-algebras, where $\exists_1, \forall_1$ are the quantifiers associated with $\langle{\bf D}_1, \s \rangle$ and $\exists_2, \forall_2$ are the quantifiers associated with $\langle {\bf D}_2, \s \rangle $.     
\end{Example}

Theorem \ref{cuantificador asociado a C} also allows us to build an example showing that the variety $\mathbb{MG}_{\s}$ is not locally finite.  Consider the $G_{\s}$-algebra $\langle [0,1]_{\mathbb{G}}, \s \rangle^{\mathbb{N}}$, where $[0,1]_{\mathbb{G}}$ is the standard Gödel algebra on the real interval and $\s$ is the unary operation 
 defined by $\s x :=1-x$. Let $C$ be the set of constant sequences in $\langle [0,1]_{\mathbb{G}}, \s \rangle^{\mathbb{N}}$. It is easy to check that  $\langle {\bf C},  \s  \rangle$ is an $m$-relatively complete subalgebra of $\langle [0,1]_{\mathbb{G}}, \s \rangle^{\mathbb{N}}$. Thus, by the above theorem we obtain that $\langle \langle [0,1]_{\mathbb{G}}, \s \rangle^{\mathbb{N}}, \exists, \forall \rangle$ is an $MG_{\s}$-algebra.   The same argument as in \cite{Complet} shows that the  subalgebra of $\langle \langle [0,1]_{\mathbb{G}}, \s \rangle^{\mathbb{N}}, \exists, \forall \rangle$ generated by the element $a(n):=1- \frac{1}{n}$, $n \in \mathbb{N}$, is infinite. This shows that  the variety $\mathbb{MG}_{\s}$ is not locally finite.

We now use Theorem \ref{cuantificador asociado a C} to prove the FEP.

\begin{Lemma} \label{inmersion finita}
If  $\langle {\bf A}, \s, \exists, \forall \rangle$ is a subdirectly irreducible ${MG}_{\s}$-algebra and $N:=\{e_1, e_2, \dots, e_n\}$ is a finite set contained in $A$, then there exists a finite subalgebra $\langle {\bf A}_1, \s \rangle$ of $\langle {\bf A}, \s \rangle$ and operations $\forall_1$ and $\exists_1$ with the following properties:

\begin{enumerate}[$(1)$]
\item $\langle {\bf A}_1, \s, \exists_1, \forall_1 \rangle$ is a finite subdirectly irreducible   ${MG}_{\s}$-algebra.

\item $N \subseteq A_1$.

\item If $a \in A_1$ and $\forall a \in A_1$, then $\forall_1 a = \forall a$.
\item If $a \in A_1$ and $\exists a \in A_1$, then $\exists_1 a = \exists a$.

\end{enumerate}
\end{Lemma}

\begin{Proof}
Let $\langle {\bf A}_1, \s \rangle$ be the subalgebra  of $\langle {\bf A}, \s \rangle$ generated by $N$. Of course, $N \subseteq A_1$ and since the variety $\mathbb{G}_{\s}$ is locally finite, $A_1$ is finite.

Let us consider  $C:= A_1 \cap \forall A$. Observe that $\langle {\bf C}, \s \rangle$ is a totally ordered subalgebra of $\langle {\bf A}_1, \s \rangle$. We claim that it is an $m$-relatively  complete subalgebra. In fact, it is clear that $(s_1)$ holds since $C$ is a finite chain. To show $(s'_2)$, suppose that  $1=c \vee a$, with $c \in C$ and $a \in A_1$. Since $\langle \exists {\bf A}, \s \rangle$ is an $m$-relatively complete subalgebra of $\langle {\bf 
A}, \s \rangle$, $c=1$ or $a=1$.

Since $\langle {\bf C}, \s \rangle$ is an $m$-relatively complete subalgebra of $\langle {\bf A}_1, \s \rangle$, by Theorem \ref{cuantificador asociado a C} we have that $\langle {\bf A}_1, \s, \exists_1, \forall_1 \rangle$ is an $MG_{\s}$-algebra where $\exists_1$ and $\forall_1$ are defined by:

$$\exists_1 a:= \min\{c \in C: c\geq a\}, \ \ \ \ \forall_1 a:= \max\{c \in C: c\leq a\},$$
Moreover, since $\langle {\bf C}, \s \rangle$ is totally ordered, $\langle {\bf A}_1, \s, \exists_1, \forall_1 \rangle$   is subdirectly irreudible by Corollary \ref{subdirectamente irreducibles MGalgebras simétricas}. 

To prove $(3)$, suppose $a \in A_1$ and $\forall a \in A_1$. Let us prove that $\forall a =\forall_1 a =\max\{c \in C: c \leq a\}$. We know that $\forall a \in C$ and $\forall a \leq a$. Suppose that $m \in C$ and $m \leq a$. Since $m \in C=\forall A \cap A_1$, there exists $m' \in A$ such that $m=\forall m'$. Then $m \leq a$ implies $\forall \forall m' \leq \forall a$ and thus $m \leq \forall a$. This shows that $\forall_1 a = \forall a$.

To prove $(4)$, let $a \in A_1$ and $\exists a \in A_1$. Since $\langle {\bf A}_1, \s \rangle$ satisfies $(Q)$ we obtain $\s \forall \s a \in A_1$. Then, $\s a, \forall \s a \in A_1$. Using $(3)$ we obtain $\forall \s a= \forall_1 \s a$ and consequently, $\exists a = \s \forall \s a = \s \forall_1 \s a = \exists_1 a$ and the proof is complete.

\end{Proof}

\begin{Theorem} \label{FEP}
$\mathbb{MG}_\sim$ has the FEP.
\end{Theorem}

\begin{Proof}
By Lemma \ref{inmersion finita} we have that the class of finite subdirectly irreducible algebras in $\mathbb{MG}_\sim$ has the F.E.P. and so $\mathbb{MG}_\sim$ has the FEP too.
\end{Proof}

\begin{Corollary} \label{MG generada por m finitos}
The variety  $\mathbb{MG}_{\s}$ is generated as a quasivarity by its finite subdirectly irreducible members.\end{Corollary}

\begin{Proof}
Let ${\cal S}_{fin}$ be the class of all
finite  subdirectly irreducible algebras in $\mathbb{MG}_{\s}$. Suppose we have  terms $t_1, \ldots, t_n, t$ in variables $x_1, \ldots, x_k$ such that ${\cal S}_{fin} \models (t_1 \approx 1 \wedge \ldots \wedge t_n \approx 1) \Rightarrow t \approx 1$ but $\mathbb{MG}_{\s} \not\models (t_1 \approx 1 \wedge \ldots \wedge t_n \approx 1) \Rightarrow t \approx 1$. Then, there exists a subdirectly irreducible algebra $\langle {\bf A}, \s, \exists, \forall \rangle \in \mathbb{MG}_{\s}$  with $\langle {\bf A}, \s, \exists, \forall \rangle \not \models (t_1 \approx 1 \wedge \ldots \wedge t_n \approx 1) \Rightarrow t \approx 1$. Let $T $ be the set of all subterms of $t_1, \ldots, t_n, t$, and let $\overline{a} := (a_1, \ldots, a_k) \in  A^k$  be such that $t_1(\overline{a}) = 1, \ldots, t_n(\overline{a}) = 1$ but $ t(\overline{a}) \ne 1$. Further,  let $\langle {\bf A}_1, \s, \exists,  \forall \rangle$ be the finite subdirectly irreducible algebra in $\mathbb{MG}_{\s}$ of Lemma \ref{inmersion finita} considering $N:=\{s(\overline{a}) : s \in T\}$. Then, we have that $\langle {\bf A}_1, \s, \exists, \forall \rangle \not\models (t_1 \approx 1 \wedge \ldots \wedge t_n \approx 1) \Rightarrow t \approx 1$, contradicting our assumption.
\end{Proof}

\section{The subvariety  $\mathbb{CMG}_{\s}$}

In this section we focus our attention on a particular subvariety of $\mathbb{MG}_{\s}$, namely $\mathbb{CMG}_{\s}$, whose members are characterized within the variety $\mathbb{MG}_{\s}$ by the following condition:

\begin{equation}
\tag{$C$}
\label{ec CMGsim algebras}
\exists (x \wedge \s x) \leq \forall (x \vee \s x). 
\end{equation}
$\mathbb{CMG}_{\s}$ is a proper subvariety of $\mathbb{MG}_{\s}$ since the subalgebras in Example \ref{ejemplos de construcción de MG algebras simetricas} with $D_2=\{0,1\}$ do not satisfy $(C)$ for $n \geq 3$.

In Section \ref{completitud} we show that the variety $\mathbb{CMG}_\sim$ is the equivalent algebraic semantics of the one-variable monadic fragment of the first-order logic ${\cal G}\forall_\sim$.  Condition $(C)$ arises from the fact that  $(\exists x) (\varphi \wedge \s \varphi) \to (\forall x) (\varphi \vee \s \varphi)$ is a theorem of ${\cal G}\forall_\sim$. To show this, take  a formula $\varphi$ in ${\cal G}\forall_\sim$, a variable $w$ such that $\varphi$ does not contain $w$  and $\psi:= \varphi[w/x]$ (the formula that results from substituting all the free occurrences of $x$ with $w$). Since ${\cal G}_\sim$ proves $(\psi \wedge \s \psi) \to (\varphi \vee \s \varphi)$ (easy consequence of \cite[Theorem 7]{EG})  we have that ${\cal G}\forall_\sim\vdash (\psi \wedge \s \psi) \to (\varphi \vee \s \varphi)$. Using generalization we obtain ${\cal G}\forall_\sim \vdash (\forall x) ((\psi \wedge \s \psi) \to (\varphi \vee \s \varphi))$ and from the fact that $\psi$ does not contain $x$ freely and $(\forall 2)$ we have that ${\cal G}\forall_\sim \vdash  (\psi \wedge \s \psi) \to (\forall x) (\varphi \vee \s \varphi)$. Now, generalization, the fact that $\varphi$ does not contain $w$ and $(\exists 2)$ yield $ {\cal G}\forall_\sim \vdash  (\exists w)(\psi \wedge \s \psi) \to (\forall x) (\varphi \vee \s \varphi)$. Finally, ${\cal G}\forall_\sim \vdash  (\exists x)(\varphi \wedge \s \varphi) \to (\forall x) (\varphi \vee \s \varphi)$  follows from  \cite[Theorem 5.1.17.]{Hj} and Modus Ponens. 
It is worth noting that even when $x$ is the only variable in $\varphi$, the proof of $(\exists x) (\varphi \wedge \s \varphi) \to (\forall x) (\varphi \vee \s \varphi)$ requires an auxiliary variable  different from $x$. Moreover, if a proof including only $x$ existed, $(C)$ should hold on the whole variety $\mathbb{MG}_{\s}$.

Next, we characterize the finite subdirectly irreducible algebras in  $\mathbb{CMG}_{\s}$ by means of the range of their monadic operators. We start by showing some basic results. Here, and subsequently, $\langle {\bf C}_{n}, \s \rangle$  denotes the $G_{\s}$-chain with $n$ elements.

\medskip

\begin{Proposition} \label{elemnto del rango entre x  y no x}
Let  $\langle {\bf A}, \s, \exists, \forall \rangle$ be a $CMG_{\s}$-algebra and $a\in A$. Then,  $a\leq \s a$ if and only if $a\leq \exists a\leq \s a$.
\end{Proposition}
\begin{Proof}
 Let $a \in A$ such that $ a \leq \s a$. Using conditions (\ref{ec CMGsim algebras}), $(Q)$ and Lemma \ref{Propiedades algebraicas de las monádicas} (2) we obtain the following equivalence: 
 $$\exists (a\land \s a)\leq  \forall (a\vee \s a) \mbox{ if and only if } a \leq \exists a \leq  \forall \s a \leq \s a.$$
\end{Proof}

An element $d$ in an $MG_{\s}$-algebra is called a {\it fixed point} if it satisfies $d = \s d$.

\begin{Corollary} \label{existe de un punto fijo es d}
If $\langle {\bf A}, \s, \exists, \forall \rangle$ is a $CMG_{\s}$-algebra with a fixed point $d$, then $\exists d=d$. 
\end{Corollary}

 Remark \ref{las si finitas son producto directo de cadenas} states that, if $\langle {\bf A}, \s, \exists, \forall \rangle$ is a finite subdirectly irreducible algebra in $\mathbb{MG}_{\s}$, then $\langle {\bf A}, \s \rangle \cong \prod_{i=1}^r \langle {\bf C}_{n_i}, \s \rangle$, for some natural numbers $n_1, \ldots, n_r$. If $\langle {\bf A}, \s,  \exists, \forall \rangle$ belongs to $\mathbb{CMG}_{\s}$, the following stronger result holds.
    
\begin{Proposition} \label{s.i. implica cadenas pares o impares}
If $\langle {\bf A}, \s, \exists, \forall \rangle $ is a finite subdirectly irreducible $CMG_{\s}$-algebra such that $\langle {\bf A}, \s \rangle =\prod_{i=1}^r \langle {\bf C}_{n_i}, \s \rangle$, then   $n_i$ is odd for all $i=1, \ldots, r$ or   $n_i$ is even for all $i=1, \ldots, r$.
\end{Proposition}

\begin{Proof}
 Assume  $\langle {\bf A}, \s, \exists, \forall \rangle$ as in the statement of the lemma, and suppose that there exists a  direct factor of $\langle {\bf A}, \s \rangle$ with an even number of elements and a direct factor with an odd number of elements. 

Observe that if $n_i$ is even, there exists an element $a_i$ in $ C_{n_i}$ such that  $\s a_i $ covers $a_i$, that is, if  $a_i \leq c \leq \s a_i$, then $c=a_i$ or $c= \s a_i$. If $n_i$ is odd, there exists an element $d \in C_{n_i}$ such that $d_i=\s d_i$.

Taking this into account, consider the element $t \in   A$ such that 

$$t(i):=\left\{ \begin{array}{lcl}
    a_i & \mbox{ if } & n_i \mbox{ is even, } \\
    d_i & \mbox{ if } & n_i \mbox{ is odd,}
  \end{array}\right.$$
\noindent where $a_i$ and $d_i$ are the elements described above.
We claim that $\exists t =t$. In fact, since by construction $t \leq \s t$, by Proposition \ref{elemnto del rango entre x  y no x} we have that $t \leq \exists t \leq \s t$. If $n_i$ is odd, $t(i) = d_i = \s d_i = \s t(i)$ and so $t(i)=(\exists t)(i)$. Suppose that $n_i$ is even, from the fact that $t \leq \exists t \leq \s t$, $(C)$ and $(Q)$ we obtain $$a_i=t(i) \leq (\exists t)(i)= (\exists (t \wedge \s t))(i) \leq \forall (t \vee \s t))(i)= (\s \exists t)(i)  \leq (\s t)(i)=\s a_i.$$ Since $\langle {\bf C}_{n_i}, \s \rangle$ has no fixed points it follows that $(\exists t)(i) < \s (\exists t)(i)$ and so $a_i \leq (\exists t)(i) < \s (\exists t)(i) \leq \s a_i$. Finally, since $\s a_i$ covers $a_i$, we have that  $t(i)=(\exists t)(i)$. Both cases show that $t(i)=(\exists t)(i)$, for all $i=1, \ldots, r$ and consequently, $\exists t =t$. Hence, $t \in \exists A$.

Consider now the elements $\s t \to t$ and $\s (\s t \to t)$ in $\exists  A$. Note that 

$$(\s t \to t)(i)= \left\{ \begin{array}{ll}
    a_i & \mbox{ if }  n_i \mbox{ is even,} \\
    1 & \mbox{ if }  n_i \mbox{ is odd,}
  \end{array}\right.  \mbox{ \  and \ } (\s(\s t \to t))(i)= \left\{ \begin{array}{ll}
    \s a_i &  \mbox{ if }  n_i  \mbox{ is even,} \\
    0 & \mbox{ if }  n_i \mbox{ is odd.}
  \end{array}\right. $$   
  
Clearly,  $\s t \to t$ and $\s (\s t \to t)$  are non-comparable elements in $\langle \exists {\bf A}, \s \rangle$ contradicting the fact that $\langle \exists {\bf A}, \s \rangle$ is totally ordered.
\end{Proof}

\begin{Corollary} \label{punto fijo cadenas impares no punto fijo cadenas pares}
 Let  $\langle {\bf A},  \s, \exists, \forall \rangle$ be  a finite subdirectly irreducible algebra in $\mathbb{CMG}_{\s}$ such that $\langle {\bf A}, \s \rangle =\prod_{i=1}^r \langle {\bf C}_{n_i}, \s \rangle$.
\begin{enumerate}[$(1)$]
\item If $\langle {\bf A}, \s, \exists, \forall \rangle$ has a fixed point, then $n_i$ is odd, for all $i=1, \ldots, r$.
\item If $\langle {\bf A}, \s, \exists, \forall \rangle$ has no  fixed points, then $n_i$ is even, for all $i=1, \ldots, r$.
\end{enumerate}
\end{Corollary}

The following result characterizes the finite subdirectly irreducible algebras $\langle {\bf A}, \s, \exists, \forall \rangle$ in the variety $\mathbb{CMG}_{\s}$ by means of $\exists  A$.

\begin{Theorem}\label{punto fijo}
  Let  $\langle {\bf A}, \s, \exists, \forall \rangle$ be  a finite subdirectly irreducible algebra in $\mathbb{MG}_{\s}$ such that $\langle {\bf A}, \s \rangle =\prod_{i=1}^r \langle {\bf C}_{n_i}, \s \rangle$.
   \begin{enumerate}[$(1)$]
     \item If $\langle {\bf A}, \s, \exists, \forall \rangle$ has a fixed point $d$, then  
           $\langle {\bf A}, \s, \exists, \forall \rangle \in \mathbb{CMG}_\sim$ if and only if  $d\in \exists A.$
     \item If $\langle {\bf A}, \s, \exists, \forall \rangle$ has no  fixed points, then $\langle {\bf A}, \s, \exists, \forall \rangle \in \mathbb{CMG}_\sim$ if and only if there exists  $t \in \exists A$ such that $(\s t)(i)$ covers $t(i)$, for all $i=1, \ldots r$.
     
   \end{enumerate}
\end{Theorem}

\begin{Proof}
  Assume  $\langle {\bf A}, \s, \exists, \forall \rangle$ as in the statement of the theorem. The forward implication in $(1)$ follows immediately from Corollary \ref{existe de un punto fijo es d}.

 Conversely, suppose that $d\in \exists A$ and  let $x\in A$. Since $d$ is a fixed point of $\langle {\bf A}, \s, \exists, \forall \rangle$, $d(i)$ is a fixed point of $\langle {\bf C}_{n_i}, \s \rangle$, for all $i=1, \ldots, r$. Thus, the coordinates of $x$ satisfy  $ x(i)\wedge \s x(i) \leq d(i)$ for all $1\leq i\leq r$. From this, Lemma \ref{Propiedades algebraicas de las monádicas}  and $(Q)$, we obtain   $\exists (x\wedge \s x)\leq \exists d =d = \s d= \s \exists d =\forall (\s d) \leq \forall (x \vee \s x).$ This proves that $\langle {\bf A}, \s, \exists, \forall \rangle$ satisfies (\ref{ec CMGsim algebras}).

  In order to prove $(2)$, suppose that  $\langle {\bf A}, \s,  \exists, \forall \rangle$ satisfies (\ref{ec CMGsim algebras}). Since $\langle {\bf A}, \s, \exists, \forall \rangle$ has no  fixed points and it is subdirectly irreducible, by Corollary \ref{punto fijo cadenas impares no punto fijo cadenas pares} we have that each factor $\langle {\bf C}_{n_i}, \s \rangle$ is a chain with an even number of elements. Then, for each $i=1, \ldots r$, there exists  $a_i \in C_{n_i}$ such that $\s a_i$ covers $a_i$.    
Let us consider the element $t$ in $A$ given by $t(i):=a_i$ for $i=1, \ldots, r$. Note that we have already proved that $t \in \exists A$ in the proof of Proposition \ref{s.i. implica cadenas pares o impares}.
  Conversely, suppose that there exists $t \in \exists A$ such that $(\s t)(i) $ covers $t(i)$ for all $i=1, \ldots, r$ and take  $x\in A$. Thus, in each factor $\langle {\bf C}_{n_i}, \s \rangle$, we have that  $x(i)\leq t(i) <  \s t(i) \leq \s x(i)$ or $\s x(i)\leq t(i)< \s t(i) \leq x(i)$. In both cases $x(i) \wedge \s x(i) \leq t(i)$, and so, $x\wedge \s x\leq t$.   Finally, since $t = \exists t$ we have that $$\exists (x\wedge \s x)\leq \exists t =t \leq  \s t= \s \exists t = \forall (\s t) \leq \forall (x\vee \s x).$$ This shows that $\langle {\bf A}, \s, \exists, \forall \rangle$ satisfies (\ref{ec CMGsim algebras}).
\end{Proof}

Our next goal is to prove that the variety $\mathbb{CMG}_{\s}$ has the FEP, and thus, is generated by its finite members. Moreover, we show that it is generated by its finite subdirectly irreducible members with a fixed point. We need to prove first  the following result for algebras in $\mathbb{MG}_\sim$.

\begin{Lemma}\label{propiedades de elementos del rango en una s.i.}
Let $\langle {\bf A}, \s, \exists, \forall \rangle$  be a finite subdirectly irreducible algebra in $\mathbb{MG}_{\s}$, where $\langle {\bf A}, \s \rangle= \prod_{i=1}^r \langle {\bf C}_{n_i}, \s \rangle $. Then:

\begin{enumerate}[$(1)$]

\item If $c \in \exists A$ and $c(i)=1$ for some $i=1, \ldots r$, then $c=1$.

\item If $c \in \exists A$ and $c(i)=\s c(i)$ for some $i=1, \ldots r$, then $c=\s c$.

\end{enumerate}

\end{Lemma}

\begin{Proof}
In order to prove (1) suppose that $c \in \exists A$ and  there exist $i,j$ such that $c(i)=1 $ and $c(j) \ne 1$. If we consider an element $x \in A$ such that $x(k)=1$ for all $k \ne i$ and $x(i)=0$, we have that $x \vee c=1$, $c \ne 1$ and $x \ne 1$. This implies that $\langle \exists {\bf A}, \s \rangle$ is not an $m$-relatively complete subalgebra of $\langle {\bf A}, \s \rangle$, which is a contradiction.

Suppose now that there exist $c \in \exists A$ and $i,j$ such that $c(i)=\s c(i)$ and $c(j) \ne \s c(j)$. It is easy to check that either $\s c \to c $ and $\s (\s c \to c)$ are non-comparable elements, or $ c \to  \s c $ and $\s (c \to \s c)$ are, which is a contradiction. Thus, $(2)$ holds.  
\end{Proof}

\begin{Lemma} \label{inmersión en punto fijo}  
\begin{enumerate}

\item[$(1)$] If $\langle {\bf A},  \s \rangle$ is a finite algebra in $\mathbb{G}_{\s}$, then $\langle {\bf A},  \s \rangle$ is embeddable in a finite algebra $\langle {\bf B}, \s \rangle$ in $\mathbb{G}_{\s}$ with a fixed point.

\item[$(2)$] If $\langle {\bf A}, \s, \exists, \forall \rangle$ is a finite subdirectly irreducible algebra in $\mathbb{MG}_{\s}$, then there exists a finite subdirectly irreducible algebra $\langle {\bf B}, \s, \exists, \forall \rangle$ in $\mathbb{CMG}_{\s}$ with a fixed point  such that  $\langle {\bf A},  \s \rangle$ is a subalgebra of $\langle {\bf B}, \s \rangle$ and $\exists A \subseteq \exists B$.
In addition, if $\langle {\bf A}, \s, \exists, \forall \rangle$ satisfies $(C)$, then     $\langle {\bf A}, \s, \exists, \forall \rangle$ is a subalgebra of $\langle {\bf B}, \s, \exists, \forall \rangle$.
 \end{enumerate}
\end{Lemma}

\begin{Proof}
To prove $(1)$, let  $\langle {\bf A}, \s \rangle$ be a finite algebra in $\mathbb{G}_{\s}$, by Lemma \ref{propiedades de Godel simetricas} (2) we can assume  that $\langle {\bf A}, \s \rangle = \prod_{i=1}^r \langle {\bf C}_{n_i}, \s \rangle$, where $\langle {\bf C}_{n_i}, \s  \rangle$ are chains in $\mathbb{G}_{\s} $. 

For each even $n_i$, let  $C_{n_i}:= \{a_1, \ldots, a_{\frac{n_i}{2}}, \s a_{\frac{n_i}{2}}, \ldots, \s a_1 \}$, where $a_1 < \ldots < a_{\frac{n_i}{2}} < \s a_{\frac{n_i}{2}} < \ldots < \s a_1$. Let us  consider the algebra $\langle \overline{{\bf C}}_{n_i}, \s \rangle$ in $\mathbb{G}
_{\s}$, where $\overline{C}_{n_i}:=C_{n_i} \cup \{d_i\}$ and  $$a_1 < \ldots < a_{\frac{n_i}{2}} < d_i = \s d_i < \s a_{\frac{n_i}{2}} < \ldots < \s a_1.$$ For each odd $n_i$ let $\langle \overline{{\bf C}}_{n_i}, \s \rangle := \langle {\bf C}_{n_i}, \s \rangle$ with fixed point $d_i$.  
It is clear that $\langle {\bf C}_{n_i}, \s \rangle$ is a subalgebra of $\langle \overline{{\bf C}}_{n_i}, \s \rangle$, and consequently $\langle {\bf A}, \s \rangle$ is a subalgebra of $\prod_{i=1}^r \langle \overline{{\bf C}}_{n_i}, \s \rangle$.  Moreover,   $(d_1, \ldots, d_r)$ is the fixed point of $\prod_{i=1}^r \langle \overline{{\bf C}}_{n_i}, \s \rangle$.

In order to prove $(2)$, let $\langle {\bf A}, \s, \exists, \forall \rangle$ be a finite subdirectly irreducible algebra in $\mathbb{MG}_{\s}$ with no  fixed points. Again, we assume that  $\langle {\bf A}, \s \rangle = \prod_{i=1}^r \langle {\bf C}_{n_i},  \s \rangle$, where each direct factor $\langle {\bf C}_{n_i},  \s \rangle$ is a chain. We proceed as the proof of $(1)$ and produce $\langle {\bf B}, \s \rangle := \prod_{i=1}^r \langle \overline{{\bf C}}_{n_i}, \s \rangle$ such that $\langle {\bf A}, \s \rangle$ is a subalgebra of $\langle {\bf B}, \s \rangle$. Here,  $d:=(d_1, \ldots, d_r)$ is the fixed point of $\langle {\bf B}, \s \rangle$ and $d \not \in A$.   

Let us prove that $ \exists A \cup \{d\}$ is the universe of an $m$-relatively complete subalgebra of $\langle {\bf B}, \s \rangle$. It is easy to prove that $ \exists  A \cup \{d\}$ is a chain from the fact that $\langle \exists \mathbf{A}, \s\rangle$ is a totally ordered subalgebra of $\langle \mathbf{A},\s\rangle$. Moreover, $\s a \in \exists A \cup \{d \}$ whenever $a \in \exists A \cup \{d\}$.  Let us prove that $\exists A \cup \{d\}$ is closed under $\to$. Using that $\langle {\bf A}, \s \rangle$ has no  fixed points and  Lemma \ref{propiedades de elementos del rango en una s.i.} (2), we obtain that if $c \in \exists A$, $c(i) \ne d_i$, for all $i=1, \ldots, r$. So, $c \leq d$ implies $c(i) < d_i$, for all $i=1, \ldots r$, and $d \leq c$ implies $d_i < c(i)$, for all $i=1, \ldots, r$. From this, we have that $c \to d, d  \to c \in \{1,d,c\}$. This shows that $ \exists A \cup \{d\}$ is a subuniverse of  $\langle {\bf B}, \s \rangle$. 

The conditions $(s_1)$ and $(s_3)$ follow immediately. Suppose now that $1= c \vee b$ with $b \in B$ and $c \in \exists A \cup \{d\}$. If  $c(i)=1$ for some $i=1, \ldots, r$, then $c \ne d$ and we have that $c \in \exists A$. Then, using Lemma \ref{propiedades de elementos del rango en una s.i.} $(1)$ we obtain that $c=1$. If  $c(i) \ne 1$ for each $i=1, \ldots, r$, since ${\bf C}_{n_i}$ is a totally ordered algebra, it must be that $b(i)=1$, for all  $i=1, \ldots, r$ and so $b=1$. This shows that $(s'_2)$ holds. We have thus proved that $ \exists  A \cup \{d\}$ is the universe of an $m$-relatively complete subalgebra of $\langle {\bf B}, \s \rangle$.

From Theorem \ref{cuantificador asociado a C} we have that $\langle {\bf B}, \s, \exists_1, \forall_1 \rangle$ is an $MG_{\s}$-algebra where $\exists_1$ and $\forall_1$ are given by
  $$\exists_1 a:=\min \{c\in \exists A\cup \{d\}: c\geq a\},$$ $$\forall_1 a:=\max \{c\in \exists A\cup \{d\}: c\leq a\}.$$ Hence, $\exists_1 B = \exists A \cup \{d\}$ and so, $\exists A \subseteq \exists_1 B$. Moreover, since $\langle \exists {\bf B}, \s \rangle$ is totally ordered and $\exists_1 d =d$, by Corollary \ref{subdirectamente irreducibles MGalgebras simétricas} and Theorem \ref{punto fijo} $(1)$ we obtain that $\langle {\bf B},\s, \exists_1, \forall_1 \rangle$ is a subdirectly irreducible $CMG_{\s}$-algebra.

Now Suppose that $\langle {\bf A}, \s, \exists, \forall \rangle$ satisfies the condition $(C)$.   
We know that  $\langle {\bf A}, \s \rangle$ is a subalgebra of $\langle {\bf B},  \s \rangle$. To prove that $\langle {\bf A}, \s, \exists,\forall \rangle$ is a subalgebra of $\langle {\bf B}, \s, \exists_1, \forall_1 \rangle$, we will show that $\forall a = \forall_1 a$ and $\exists a = \exists_1 a$, whenever $a \in A$.

Let   $a \in A$. If $a(i) \leq d_i$, for all $i$, we have that $a \leq d$ and so, $a \leq d = \s d \leq \s a$. Using that $\langle {\bf A}, \s,  \exists, \forall \rangle$ satisfies $(C)$, we obtain $\exists a \leq \forall (\s a) = \s \exists a$, and the fact that $\exists  A \cup \{d\} $ is a chain implies $\exists a \leq d$. Hence,  $d > \exists a \geq a$, and so  $\exists_1 a = \exists a$. If $a(i) \not \leq d_i$, for some $i$, we have that, $ a \not \leq d$, and consequently, $\exists_1 a = \exists a$. 
Taking this into account and using the identity $(Q)$, it follows immediately that $\forall_1 a = \forall a$, for each $a \in A$.
\end{Proof}

\begin{Proposition} \label{inmersion finita2}
If  $\langle {\bf A}, \s, \exists, \forall \rangle$ is a subdirectly irreducible $CMG_{\s}$-algebra and  $N:=\{e_1, e_2, \dots, e_n\}$ is a finite set contained in $A$, then there exists a finite subdirectly irreducible  $CMG_{\s}$-algebra $\langle {\bf A}_2, \s_2, \exists_2, \forall_2 \rangle$ with the following properties:
\begin{enumerate}[$(1)$]
\item The subalgebra of $\langle {\bf A}, \s \rangle $ generated by $N$ is a subalgebra of  $\langle  {\bf A}_2, \s_2 \rangle$.

\item If $a \in N$ and $\exists a \in A_2$, then $\exists_2 a = \exists a$.

\item If $a \in N$ and $\forall a \in A_2$, then $\forall_2 a = \forall a$.

\end{enumerate}
\end{Proposition}

\begin{Proof}
 Let $\langle {\bf A}, \s, \exists, \forall \rangle$ be a subdirectly irreducible $CMG_{\s}$-algebra. By Lemma \ref{inmersion finita} there exists a finite subdirectly irreducible algebra $\langle {\bf A}_1, \s,  \exists_1, \forall_1 \rangle$ in $\mathbb{MG}_{\s}$ that satisfies the following:
 
 \begin{itemize}  
 \item $N\subseteq A_1$, 
 \item $\langle {\bf A}_1, \s \rangle$ is a subalgebra of $\langle {\bf A}, \s \rangle$, 
\item if $a \in A_1$ and $\forall a \in A_1$, then $\forall_1 a = \forall a$,
\item if $a \in A_1$ and $\exists a \in A_1$, then $\exists_1 a = \exists a$.
 
 \end{itemize}
We can further assume that $\langle {\bf A}_1, \s \rangle=\prod_{i=1}^{r} \langle {\bf C}_{n_i}, \s\rangle$.

Suppose first that $\langle {\bf A}_1, \s, \exists_1, \forall_1 \rangle$ has a fixed point $d$. Then, $d$ is also a fixed point of $\langle \mathbf{A}, \s, \exists,\forall\rangle$. Since $\langle {\bf A}, \s, \exists, \forall \rangle \in \mathbb{CMG}_{\s}$, by  Corollary \ref{existe de un punto fijo es d} we have that $\exists d= d$, and so $\exists_1 d = \exists d = d$.   Then, by Theorem \ref{punto fijo} (1), $\langle {\bf A}_1, \s, \exists_1, \forall_1 \rangle$ is a $CMG_{\s}$-algebra, and we are finished.
 
 Now suppose $\langle \mathbf{A}_1, \s, \exists_1,\forall_1\rangle$ has no fixed points. By Lemma  \ref{inmersión en punto fijo} (2), there exists a finite subdirectly irreducible  $CMG_{\s}$-algebra $\langle {\bf A}_2, \s, \exists_2, \forall_2 \rangle$ with a fixed point $d:=(d_1, \ldots, d_r)$,  such that $\langle {\bf A}_1, \s \rangle$ is a subalgebra of  $\langle {\bf A}_2, \s \rangle=\prod_{i=1}^r \langle \overline{{\bf C}}_{n_i}, \s \rangle$,  each $\langle \overline{{\bf C}}_{n_i}, \s \rangle$ is built as in the proof of Lemma \ref{inmersión en punto fijo} (1),  $\exists_2 A_2 := \exists_1 A_1 \cup \{d\}$ and $\exists_2,\forall_2$ are given by
  $$\exists_2 a:=\min \{c\in \exists A_1\cup \{d\}: c\geq a\},$$ $$\forall_2 a:=\max \{c\in \exists A_1\cup \{d\}: c\leq a\}.$$
  
 Let us prove $(2)$. Let $a \in N$ and $\exists a \in A_2$. Clearly, $\exists a \in A_1$ and so, $\exists_1 a = \exists a$. If $a(i) \leq d_i$, for all $i$, we have that $a \leq d$ and so, $a \leq d = \s d \leq \s a$. Using that $\langle {\bf A}, \s, \exists, \forall \rangle$ satisfies $(C)$, we obtain $\exists_1 a =\exists a \leq \forall (\s a) = \s \exists a= \s \exists_1 a$, and the fact that $\exists  A_1 \cup \{d\} $ is a chain implies $\exists_1 a \leq d$. Hence,  $d > \exists_1 a \geq a$, and so  $\exists_2 a = \exists_1 a= \exists a$. If $a(i) \not \leq d_i$, for some $i$, we have that, $ a \not \leq d$, and consequently, $\exists_2 a = \exists_1 a= \exists a$. 
Taking this into account and using the identity $(Q)$, it follows immediately that $\forall_2 a =\forall_1 a = \forall a$, for each $a \in N$ and $\exists a \in A_2$.
\end{Proof}

\begin{Theorem} \label{FEP2}
$\mathbb{CMG}_\sim$ has the FEP.
\end{Theorem}

\begin{Corollary}
The variety $\mathbb{CMG}_{\s}$ is generated as a quasivariety by its finite subdirectly irreducible members.
\end{Corollary}


This corollary  and Lemma \ref{inmersión en punto fijo} (2) yield the following  stronger result.

\begin{Corollary} \label{CMG generada como quasiv por si con punto fijo}
The variety $\mathbb{CMG}_{\s}$ is generated as a quasivariety by its finite subdirectly irreducible members with a fixed point.

\end{Corollary}

\section{Functional monadic $ G_{\s}$-algebras}

In this section we consider an important example of  $CMG_{\s}$-algebras. Our interest in these algebras arises from the fact that they encode the information of the models used to interpret the one-variable monadic fragment of the first-order logic ${\cal G}\forall_{\s}$. We call them {\it functional monadic $G_{\s}$-algebras} and we prove that they generate  the variety $\mathbb{CMG}_{\s}$.

We define functional monadic $G_{\s}$-algebras in the same way that functional monadic $BL$-algebras are defined in  $\cite{Bahia}$ and \cite{Complet}. Let $\langle {\bf A}, \s\rangle$ be a totally ordered $G_{\s}$-algebra  and let $X$ be a non-empty set. Consider the $G_{\s}$-algebra $\langle {\bf A}^X, \s \rangle:= \langle {\bf A}, \s \rangle^X$ and let $S$ be  the set of functions $f \in A^X$ such that $\mbox{inf}\{f(x): x \in X\}$ and $\mbox{sup}\{f(x): x \in X\}$ both exist in ${\bf A}$. For every $f \in S$, we define: $(\forall_\wedge f)(x) := \mbox{inf}\{f(y):y \in X \}$ and  $(\exists_\vee f)(x) := \mbox{sup}\{f(y):y \in X \}$. Note that $\forall_\wedge f$ and $\exists_\vee f$ are constant maps.  We write $C$ for the set of constant maps of $A^X$. Let $\langle {\bf B}, \s \rangle$ be a subalgebra  of $\langle {\bf A}^X, \s \rangle$ such that $B \subseteq S$ and $B$ is closed under $\exists_\vee$ and $\forall_\wedge$. In \cite{Bahia}  it was shown that  the structure $\langle {\bf B}, \exists_\vee, \forall_\wedge \rangle$ is a monadic $BL$-algebra, $ C \cap  B$ is the universe of an $m$-relatively complete subalgebra of ${\bf B}$ and $ C \cap  B = \exists_\vee B = \forall_\wedge B$. We claim that $\langle {\bf B}, \s, \exists_\vee, \forall_\wedge \rangle$ is an $MG_{\s}$-algebra. Indeed, it is clear that $\langle {\bf C} \cap {\bf B}, \s \rangle$ is a $G_{\s}$-subalgebra of ${\bf B}$ and so, it is  an $m$-relatively complete subalgebra of $\langle {\bf B}, \s \rangle$. Thus,  by Theorem \ref{cuantificador asociado a C},  $\langle {\bf B}, \s, \exists_\vee, \forall_\wedge \rangle$ is an $MG_{\s}$-algebra. We call  the algebras $\langle {\bf B}, \s, \exists_\vee, \forall_\wedge \rangle$ that arise in this way {\it functional monadic $G_{\s}$-algebras.} The main example of these algebras is obtained considering $\langle {\bf A}, \s \rangle$  to be the  standard $G_{\s}$-algebra $\langle [0,1]_{\mathbb{G}}, \s \rangle$ and $X=\mathbb{N}$.

\begin{Proposition}

Every functional monadic $G_{\s}$-algebra satisfies {\rm(}\ref{ec CMGsim algebras}{\rm)}.
\end{Proposition}
\begin{Proof}
Let  $\langle {\bf B}, \s, \exists_\vee, \forall_\wedge \rangle$ be a functional monadic $G_{\s}$-algebra whit $\langle {\bf B}, \s \rangle \leq \langle {\bf A}^X,  \s \rangle$ for a suitable totally ordered $G_\sim$-algebra $\langle {\bf A}, \s\rangle$ and a non-empty set $X$. Given $f \in B$, since by Lemma \ref{propiedades de Godel simetricas} (1) $\langle{\bf A}, \s \rangle$   satisfies the Kleene condition,  we have that  $f(x) \wedge \s f(x) \leq f(y) \vee \s f(y)$, for all $x,y \in X$. 
Therefore, we have that
$$\exists_\wedge (f \wedge \s f)= \mbox{sup}\{f(x) \wedge \s f(x), x \in X\} \leq \mbox{inf}\{f(y) \vee \s f(y), y \in X\}= \forall_\wedge (f \vee \s f).$$
\end{Proof}

Next we show that all finite subdirectly irreducible $CMG_{\s}$-algebras with a fixed point are isomorphic to  a functional monadic $ G_{\s}$-algebra. We first prove this result for all finite totally ordered $CMG_{\s}$-algebras.

\begin{Lemma} \label{cadena impar isomorfa a funcional}
Every finite totally ordered ${CMG}_{\s}$-algebra with a fixed point is isomorphic to a functional monadic $ G_{\s}$-algebra.
\end{Lemma}

\begin{Proof}
Let  $\langle {\bf A}, \s, \exists, \forall \rangle$ be a finite $CMG_{\s}$-chain with an odd number of elements and $\exists  A :=\{c_0, c_1, \ldots, c_{k}, \s c_{k-1}, \ldots, \s c_0\}$ where  $  0=c_0 < c_1 < \ldots < c_k=\s c_k < \s c_{k-1} < \ldots < \s c_0=1$.

Put $A_i:=\{x \in A:  c_{i-1} < x < c_i\}$ and $\s A_i:=\{x \in A:  \s c_i < x < \s c_{i-1}\}$  for $1 \leq i \leq k$, we have that $\exists A_i = \{c_i\}$,  $\forall A_i = \{c_{i-1}\}$, $\exists (\s A_i) = \{ \s c_{i-1}\}$ and  $\forall (\s A_i) = \{ \s c_i\}$ whenever $A_i \ne \emptyset$.

Consider the functional monadic $G_{\s}$-algebra $\langle [0,1]_{\mathbb{G}}^\mathbb{N}, \s, \exists_\vee,  \forall_\wedge \rangle$. We will prove that $\langle{\bf A}, \s, \exists, \forall \rangle$ is embeddable in $\langle [0,1]_{\mathbb{G}}^\mathbb{N}, \s,  \exists_\vee, \forall_\wedge \rangle$. 
In order to show that, take the constant sequences $g_i$ in $[0,1]_{\mathbb{G}}^\mathbb{N}$ given by $$g_i(n):=\frac{i}{2k},$$   for each $ n \in \mathbb{N}$ and  $0 \leq i \leq k$.  

Note that $g_k$ is the fixed point in $\langle [0,1]_{\mathbb{G}}^\mathbb{N}, \s \rangle$ and $0 \equiv g_0 < g_1 < \ldots < g_k \equiv \frac{1}{2}$. In addition, for $1 \leq i \leq k$ and $j\geq1$, consider the  maps $f_i^{(j)} \in [0,1]_{\mathbb{G}}^\mathbb{N}$  such that

$$f_i^{(j)}(n) := \left\{ \begin{array}{ll} \frac{i}{2k} - \frac{1}{2k} \left( 1- \frac{1}{n+1} \right)^j & \mbox{ if } n \mbox{ es even,} \\
\frac{i}{2k} - \frac{1}{2k} \left( \frac{1}{n+1} \right)^j & \mbox{ if } n \mbox{ is odd.} \\  \end{array} \right.$$

Observe that

$$0\equiv  g_0 < f_i^{(j)}  < g_{k}\equiv \frac{1}{2} \mbox{ \ \ \ and \ \ \ } g_{i-1} < f_i^{(1)} < f_i^{(2)} < \cdots < f_i^{(j)}< f_i^{(j+1)} < \cdots <  g_{i}.$$ Thus, 

 $$ \frac12 <  \s f_i^{(j)}  < 1  \mbox{ \ \ \ and \ \ \ } \s g_{i} < \cdots < \s f_i^{(j+1)} < \s f_i^{(j)} < \cdots < \s f_i^{(2)}< \s f_i^{(1)} <  \s g_{i-1}. $$

Moreover, an easy computation shows that $\forall_\wedge f_i^{(j)} = g_{i-1}$, $\forall_\wedge \s f_i^{(j)} =  \s g_{i} = 1-g_{i}$, $\exists_\vee f_i^{(j)} = g_{i}$ and $\exists_\vee \s f_i^{(j)} =  \s g_{i-1} = 1-g_{i-1}.$ 

Finally, let $A_i:= \{ a_i^{(1)}, a_i^{(2)}, \ldots, a_i^{(r_i)}\}$ and
$\s A_i:= \{ \s a_i^{(r_i)}, \ldots, \s a_i^{(1)}\}$ for $1 \leq i \leq k$ where $a_i^{(1)} < a_i^{(2)} < \ldots < a_i^{(r_i)}$ and $ \s a_i^{(r_i)} < \s a_i^{(r_i-1)} < \ldots <\s a_i^{(1)} $ and consider the map $\varphi: A \to [0,1]_{{\cal \mathbb{G}}}^\mathbb{N}$ given by

$$ \left\{ \begin{array}{ll} \varphi(c_i):=g_i & \mbox{ for } 0 \leq i \leq k, \\
\varphi(\s c_i):=\s g_i & \mbox{ for } 0 \leq i \leq k, \\
\varphi(a_i^{(j)}):=f_i^{(j)} & \mbox{ for } 1 \leq i \leq k \mbox{ and } 1 \leq j \leq r_i, \\
\varphi(\s a_i^{(j)}):=\s f_i^{(j)} & \mbox{ for } 1 \leq i \leq k \mbox{ and } 1 \leq j \leq r_i. \\
  \end{array} \right.$$

It follows  immediately that  $\varphi$ is an embedding. 
\end{Proof}

In \cite{Complet} the authors  built a monadic Gödel algebra $\langle {\bf A}^n, \overline{\exists}, \overline{\forall} \rangle$ by means of a monadic Gödel chain $\langle {\bf A}, \exists, \forall \rangle$ and a natural number $n$ where the quantifiers are defined as follows: 
 $$\overline{\exists}(a_1, \ldots, a_n)(i):= \exists (a_1 \vee \ldots \vee a_n)    \mbox{ \ \ \ and  \ \ \ } \overline{\forall}(a_1, \ldots, a_n)(i):= \forall (a_1 \wedge \ldots \wedge a_n),$$ for $1 \leq i \leq n$. Here, $\overline{\exists}(A^n) = \overline{\forall}( A^n) = \{(c, \ldots, c): c \in \exists A\}$ is the universe of an $m$-relatively complete subalgebra of ${\bf A}^n$. 
 
 If we now consider  a finite $CMG_{\s}$-chain $\langle {\bf A}, \s, \exists, \forall \rangle$   with a fixed point, we have that $ \langle {\bf A}^n,  \s \rangle= \langle {\bf A}, \s \rangle^n$  is a $G_{\s}$-algebra.  A straightforward  verification shows that $\langle \overline{\exists}{\bf A}^n, \s \rangle$ is an $m$-relatively complete subalgebra of $\langle {\bf A}^n,  \s \rangle$ that contains the fixed point of  ${\bf A}^n$. Thus, we obtain that $ \langle {\bf A}^n, \s, \overline{\exists}, \overline{\forall} \rangle$ is an algebra in the variety $\mathbb{CMG}_{\s}$.

\begin{Lemma} \label{construcción de CMG_sim algebra}
For each finite $CMG_{\s}$-chain $\langle {\bf A}, \s, \exists, \forall \rangle$ with a fixed point  and every natural number $n$, the algebra $\langle {\bf A}^n, \s, \overline{\exists}, \overline{\forall} \rangle$ is isomorphic to a functional monadic $G_{\s}$-algebra.
\end{Lemma}

\begin{Proof}
Let  $\langle {\bf A}, \s, \exists, \forall \rangle$ be a finite totally ordered $CMG_{\s}$-algebra  with a fixed point. By the previous lemma   there exists an  embedding $\varphi$ from $\langle {\bf A}, \s, \exists, \forall \rangle$ to $ \langle {\bf B}^X, \s, \exists_\vee, \forall_\wedge \rangle$ where $\langle {\bf B}, \s \rangle$ is a  $G_{\s}$-chain. Thus, $\varphi(\exists a)= \exists_\vee \varphi(a)$, $\varphi(\forall a)= \forall_\wedge \varphi(a)$ and $\varphi(\s a) = \s \varphi(a)$, for all $a \in A$.    
It was proved in \cite{Complet}  that $\varphi$ induces an embedding $\psi:  {\bf A}^n \to {\bf B}^{X \times \{1, \ldots, n\}}$ given by $$\psi(a_1, \ldots, a_n)(x,i) := \varphi(a_i)(x)$$ such that $\psi(\overline{\exists}(a_1, \ldots, a_n)) = \exists_\vee \psi(a_1, \ldots, a_n)$ and $\psi(\overline{\forall}(a_1, \ldots, a_n)) = \forall_\wedge \psi(a_1, \ldots, a_n)$. We claim that $\psi$ also satisfies $ \psi(\s (a_1, \ldots, a_n))=\s \psi(a_1, \ldots, a_n)$. Indeed, 
\begin{align*}
\s \psi(a_1, \ldots, a_n)(x,i) & = \s \varphi(a_i)(x)\\
& = \varphi(\s a_i) (x) \\
&= \psi(\s (a_1, \ldots, a_n))(x,i). 
\end{align*}
Therefore, $\langle {\bf A}, \s, \overline{\exists}, \overline{\forall} \rangle$ is isomorphic to a functional monadic $G_{\s}$-algebra.
\end{Proof}

The following basic property of finite subdirectly irreducible $CMG_\sim$-algebra is the last ingredient we need to prove the main theorem of this section.

\begin{Lemma}
Let $\langle {\bf A}, \s, \exists, \forall \rangle$ be  a  finite subdirectly irreducible $CMG_{\s}$-algebra with a fixed point where $\langle {\bf A}, \s \rangle =\prod_{i=1}^r \langle {\bf A}_i, \s \rangle $. Then,  $\pi_i|_{\exists {\bf A}}$ is injective, where $\pi_i:{\bf A} \to {\bf A}_i$ are the projection maps, for every $i \in \{1, \ldots, r\}$. Thus,  $|\exists {\bf A}| =|\pi_i({\exists {\bf A}})|$.
\end{Lemma}

\begin{Proof}
Suppose that there exist $a,b \in \exists A$ such that $a \ne b$ and $a(i)=b(i) $ for some $i \in \{1,2, \ldots, r\}$. Observe that $a(i)=b(i) \ne 1$.  Since otherwise, by Lemma \ref{propiedades de elementos del rango en una s.i.} (1) we have that $a=b=1$ which is a contradiction. Since $\exists {\bf A}$ is totally ordered we can suppose, without loss of generality, that $a \leq b$ and from $a \ne b$ we have that $a(j) <b(j)$ for some $j\in \{1, \ldots r\}$. We claim that the elements $b \to a$ and $b$ are no comparable. Indeed, we have that $b(i) < (b \to a)(i)=1$ and  $(b \to a)(j)=a(j)< b(j)$. This completes the proof.
\end{Proof}

\begin{Theorem}

Every finite subdirectly irreducible  ${CMG}_{\s}$-algebra  with a fixed point is isomorphic to a monadic functional $G_{\s}$-algebra.

\end{Theorem}

\begin{Proof}
 Let $\langle {\bf A}, \s, \exists, \forall\rangle$ be a finite subdirectly irreducible $CMG_{\s}$-algebra with a fixed point. By Proposition \ref{s.i. implica cadenas pares o impares}, we can assume  there are finite $G_{\s}$-chains ${\bf A}_{i}$ with a fixed point such that $\langle {\bf A}, \s \rangle = \prod_{i=1}^r \langle {\bf  A}_i, \s \rangle $. By the previous lemma, we can assume without loss of generality that $C:=\pi_i(\exists A)=  \{ c_0, c_1, \ldots, c_{k-1}, c_k, \s c_{k-1}, \ldots, \s c_0\}$, for all $i \in \{1, \ldots, r \}$, where  $$  0=c_0 < c_1 < \ldots < c_k=\s c_k < \s c_{k-1} < \ldots < \s c_0=1.$$ We further assume $A_i \cap  A_j =  C$, for $i \ne j$. Note that $\exists A = \{(x,x, \ldots, x) \in A: x \in   C\}$.

Put $A_{ij}:=\{x \in A_i:  c_{j-1} < x < c_j\}$ and $\s A_{ij}:=\{x \in A_i:  \s c_{j} < x < \s c_{j-1}\}$  for $1 \leq i \leq r$, $1 \leq j \leq k$.  We can write the underlying structure of distributive lattice of ${\bf A}_i$ in terms of lattice ordinal sum as follows:
$$\{c_0\} \oplus A_{i1} \oplus \{c_1\} \oplus  \ldots \oplus A_{ik} \oplus \{c_{k}\}  \oplus \s A_{ik} \oplus  \ldots \oplus \{\s c_1\} \oplus \s A_{i1} \oplus \{\s c_0\}. $$

Let  $B:=  \bigcup_{i=1}^r A_i$ and consider the Gödel algebra ${\bf B}$ with universe $B$ given by the following ordinal sum:

$$ \{c_0\} \oplus \bigoplus_{i=1}^r A_{i1} \oplus \{c_1\} \oplus  \ldots \oplus \bigoplus_{i=1}^r A_{ik} \oplus \{c_{k}\}  $$ $$
 \oplus \bigoplus_{i=0}^{r-1} \s A_{(r-i)k} \oplus  \ldots \oplus \{\s c_1\} \bigoplus_{i=0}^{r-1} \s A_{(r-i)1} \oplus \{\s c_0\}.   $$

Let $\sim$ be the natural involutive anti-isomorphism such that  $\langle {\bf B}, \s \rangle$ is a  $G_{\s}$-algebra. In addition, for $i=1, \ldots r$, we have that $\langle {\bf A}_i, \s \rangle$ is a $G_{\s}$-subalgebra of $\langle {\bf B}, \s \rangle$.

It is easy to check that $C$ is the universe of an $m$-relatively complete subalgebra of $\langle {\bf B}, \s \rangle$. Hence, by Theorems \ref{cuantificador asociado a C} and \ref{punto fijo} we can define quantifiers in ${\bf B}$ given by $\exists x := \min\{ c \in C: x \leq c\}$ and $\forall x := \max\{ c \in C: c \leq x\}$ such that  $\langle {\bf B}, \s,  \exists, \forall \rangle$ is a $CMG_{\s}$-algebra.
Therefore, we can apply the construction described before  Lemma \ref{construcción de CMG_sim algebra} and obtain the $CMG_{\s}$-algebra  $\langle {\bf B}^r, \s, \overline{\exists}, \overline{\forall} \rangle$.

Since $\langle {\bf A}_i, \s \rangle$ is a subalgebra of $\langle {\bf B}, \s \rangle$, we obtain that $\langle {\bf A}, \s \rangle= \langle \prod_{i=1}^r {\bf  A}_i, \s \rangle $ is a subalgebra of $ \langle {\bf  B}^r, \s \rangle$.  We claim that  $\langle {\bf A}, \s, \exists, \forall \rangle$ is a subalgebra of $\langle {\bf B}^r, \s, \overline{\exists}, \overline{\forall} \rangle$. Indeed, for $(a_1, \ldots, a_r )\in A$ we have 
\begin{align*} 
(\exists (a_1, \ldots, a_r ))(i) & = \min\{c \in C: a_j \leq c, \mbox{ for all } j=1, \ldots,  r\} \\
& = \min\{c \in C: a_1 \vee \ldots \vee a_r  \leq c\} \\
& = (\overline{\exists} (a_1, \ldots, a_r ))(i).
\end{align*}
In the same way we obtain that $\forall a= \overline{\forall}a$, for all $a \in A$.  Finally, since by Lemma \ref{construcción de CMG_sim algebra} we have that  $\langle {\bf B}^r, \s, \overline{\exists}, \overline{\forall} \rangle$ is isomorphic to a functional monadic $G_{\s}$-algebra, so is $\langle {\bf A}, \s, \exists, \forall \rangle $.
\end{Proof}

This theorem and Corollary \ref{CMG generada como quasiv por si con punto fijo} yield the following result.

\begin{Corollary} \label{CMG generada como cuaisvariedad por funcionales con punto fijo}
$\mathbb{CMG}_{\s}$ is generated as a quasivariety by its  finite functional monadic $G_{\s}$-algebras with a fixed point.
\end{Corollary}


Taking into account the construction given in the proof of Lemma \ref{cadena impar isomorfa a funcional} we obtain the following stronger result.

\begin{Corollary} $\mathbb{CMG}_{\s}$ is generated as a quasivariety by the functional monadic $G_{\s}$-algebra $\langle [0,1]_{\mathbb{G}}^\mathbb{N}, \s, \exists_\vee, \forall_\wedge \rangle$.
\end{Corollary}

\section{ A completeness result} \label{completitud}

  In this section we show that the class $\mathbb{CMG}_{\s}$ is the equivalent  algebraic semantics of the one-variable monadic fragment of the  predicate 
logic ${\cal G}\forall_{\s}$. 
To prove this, we introduce the logics ${\cal S}5({\mathbb G}_{\s})$, defined semantically, and the Hilbert-style syntactic calculus  
${\cal S}5({\cal G}_{\s})$.   
We show that ${\cal S}5({\cal G}_{\s})$ is complete with respect to the class $\mathbb{CMG}_{\s}$ and ${\cal S}5({\mathbb G}_{\s})$ is equivalent  to the one-variable monadic fragment of the first-order logic ${\cal G}\forall_{\s}$. Finally, using a  completeness result via functional algebras  analogous to the one given in \cite{Complet} we obtain that the equivalent  algebraic semantics of the one-variable monadic fragment of the predicate logic ${\cal G}\forall_{\s}$ is  the variety $\mathbb{CMG}_{\s}$.

Next, we define the logic ${\cal S}5({\mathbb G}_{\s})$ in the same way that P. Hájek defined the logic ${\cal S}5(\mathbb{BL})$ in \cite{Hj}.  The language of this logic is the language of  the logic ${\cal G}_{\s}$ augmented with the unary connectives $\square$ and $\lozenge$. The formulas in ${\cal S}5({\mathbb G}_{\s})$ are interpreted in a structure given by a triple ${\bf K}=\langle W, e, {\bf L}\rangle$ where $W$ is a non-empty set, ${\bf L}$ is a totally ordered $G_{\s}$-algebra  and $e: W \times Prop \to L$ is a function defined on  pairs $(w,p)$ consisting of an element $w \in W$ and a propositional variable $p$. The truth degree $||\varphi||_{{\bf K}, w}$ of a formula $\varphi$ in ${\bf K}$ at $w$ is defined by setting $||\varphi||_{{\bf K}, w}=e(w, \varphi)$ if $\varphi \in Prop$ and extending  its definition for any formula $\varphi$ as follows:

\begin{itemize}
\item $||0||_{{\bf K}, w} :=  0^{\bf L}$, $||1||_{{\bf K}, w} := 1^{\bf L}$,

\item $||\varphi \wedge \psi||_{{\bf K}, w} := ||\varphi ||_{{\bf K}, w} \wedge^{\bf L} || \psi||_{{\bf K}, w}$, and the same for $\vee$ and $\to$,

\item $||\s \varphi||_{{\bf K}, w}:= \s^{\bf L} ||\varphi ||_{{\bf K}, w},$

\item $||\square \varphi ||_{{\bf K}, w}:=\inf\{||\varphi ||_{{\bf K}, v}: v \in W\}$,

\item $||\lozenge \varphi ||_{{\bf K}, w}:=\sup\{||\varphi ||_{{\bf K}, v}: v \in W\}$.

\end{itemize}

Observe that $||\square \varphi ||_{{\bf K}, w}$ and $||\lozenge \varphi ||_{{\bf K}, w}$ may be undefined. Hence, we only consider  {\it safe} structures, that is,  structures ${\bf K}$ where  $||\square \varphi ||_{{\bf K}, w}$ and $||\lozenge \varphi ||_{{\bf K}, w}$ are defined for any formula $\varphi$.

We now define the consequence relation $\models_{{\cal S}5({\mathbb G}_{\s})}$.  Given a set of formulas $\Gamma$, a {\it model} of $\Gamma$ is a safe structure ${\bf K} =\langle W, e, {\bf L}\rangle$ such that $||\varphi ||_{{\bf K}, w}=1$ for every $\varphi \in \Gamma$   and every $w \in W$. Thus, given a set of formulas $\Gamma$ we define
$$ \Gamma \models_{{\cal S}5({\mathbb G}_{\s})} \varphi \mbox{ if and only if }  ||\varphi ||_{{\bf K}, w}=1 \mbox{ for every model } {\bf K}=\langle W, e, {\bf L} \rangle \mbox{ of }  \Gamma \  \mbox{ and every } w \in W.$$

In \cite{Hj} the author  proved that the modal logic ${\cal S}5(\mathbb{BL})$ is equivalent to the monadic fuzzy predicate calculus $m{\cal BL}\forall$ with unary predicates and just one object variable $x$ . The atomic formula $P_i(x)$ corresponds to the propositional variable $p_i$, modalities $\square$ and $\lozenge $ correspond to the  quantifiers $\forall x$ and $\exists x$, respectively, and the triples ${\bf K} =\langle W, e, {\bf L}\rangle$ correspond to models for the $m{\cal BL}\forall$.   
In the same way we can obtain that ${\cal S}5(\mathbb{G}_{\s})$ is actually equivalent to the one-variable  monadic fragment of the predicate logic ${\cal G}\forall_{\s}$  because there is a natural correspondence between formulas of both logics, between corresponding models and between corresponding truth degrees.

We now introduce the Hilbert-style syntactic calculus ${\cal S}5({\cal G}_{\s})$ as an expansion of ${\cal S}5'(BL)$ given in \cite{Bahia} enriched with a unary connective  $\s$. The axioms of ${\cal S}5({\cal G}_{\s})$ are the axioms in ${\cal S}5'(BL)$ plus the following axiom schemata:

\begin{enumerate}
\item[$(G)$] $\varphi \to \varphi \& \varphi$

\end{enumerate}

\begin{enumerate}[$(\s1)$]
  \item $\s(\s\varphi)\equiv\varphi$
  \item $\neg\varphi\to\s\varphi$
  \item $\Delta(\varphi\to\psi)\to\Delta(\s\psi\to\s\varphi)$
\end{enumerate}
\begin{enumerate}[$(\Delta1)$]
  \item $\Delta\varphi\lor\neg\Delta\varphi$
  \item $\Delta(\varphi\lor\psi)\to(\Delta\varphi\lor\Delta\psi)$
  \item[$(\Delta5$)] $\Delta(\varphi\to\psi)\to(\Delta\varphi\to\Delta\psi)$
\end{enumerate}
\begin{enumerate}[$(\lozenge1)$]
\item $\square \varphi \ \equiv \ \s \lozenge \s \varphi$
 \item $\lozenge (\varphi \wedge \s \varphi) \to \ \square (\varphi \vee \s \varphi)$
\end{enumerate}
\noindent where $\neg\varphi$ is $\varphi\to\overline0$ and $\Delta\varphi$ is $\neg\s\varphi$, and the inference rules of ${\cal S}5({\cal G}_{\s}) $  are those of ${\cal S}5'(BL)$ (modus ponens and generalization) and necessitation for $\Delta$ ($\displaystyle \frac{\varphi}{\Delta \varphi}$). Recall that $(G)$ implies that  $\varphi \& \psi$ is equivalent to $\varphi \wedge \psi$.

\begin{Proposition}
${\cal S}5({\cal G}_{\s})$ is an implicative logic.
\end{Proposition}

\begin{Proof}
Since ${\cal S}5'(BL)$ is an implicative logic (see \cite[Lemma 4.7]{Bahia}) it  only remains to prove  that the condition  $\varphi \to \psi, \psi \to \varphi \vdash \s \psi \to \s \varphi$ holds in ${\cal S}5({\cal G}_{\s})$.  In fact,  by necessitation, $\varphi \to \psi \vdash \Delta (\varphi \to \psi)$ and using $(\s 3)$ we obtain  $\varphi \to \psi \vdash \Delta (\s \psi  \to \s \varphi)$. Finally, using the definition of $\Delta$ and $(\s 2)$ we have that $\varphi \to \psi \vdash \s (\s (\s \psi \to \s \varphi))$ and so, $\varphi \to \psi \vdash \s \psi \to \s \varphi$.    
\end{Proof}

Recall that given an implicative logic ${\cal L}$ in the language $L$ it is possible to define the class $Alg^\ast{\cal L}$ of all algebras ${\bf A}$ of type $L$  with an element $1 \in A$ that satisfy:

\begin{enumerate}
\item[$(Alg1)$] For all $\Gamma \subseteq Fm(L)$ and $h$ homomorphism from $Fm(L)$ to ${\bf A}$, if $\Gamma \vdash_{\cal L} \varphi$ and $h[\Gamma] \subseteq \{1\}$ then $h(\varphi)=1$. 

\item[$(Alg2)$] For all $a,b \in A$, if $a \to b=1$ and $b \to a =1$, then $a=b$.

\end{enumerate}
Moreover, it is well known that  $Alg^\ast{\cal L}$ is the equivalent algebraic semantics of the logic ${\cal L}$. Hence, we have that

\begin{eqnarray} \label{semantica algebraica equivalente}
  \Gamma \vdash_{{\cal L}} \varphi
 \mbox{ if and only if }  \Gamma \models_{Alg^\ast{\cal L}} \varphi,
\end{eqnarray}

\noindent where $\Gamma \models_{Alg^\ast{\cal L}} \varphi$  means that $h(\varphi)=1$, for every  homomorphism  $h$ from $Fm(L)$ to ${\bf A}$ such that $h[\Gamma] \subseteq \{1\}$.

The following result  characterizes the class $Alg^\ast{\cal S}5({\cal G}_{\s})$.

\begin{Proposition} \label{completitud fuerte con S5 y CMG}
The equivalent algebraic semantics of ${\cal S}5({\cal G}_{\s})$ is the variety $\mathbb{CMG}_{\s}$.
\end{Proposition}

\begin{Proof}
In \cite{EG}, F. Esteva, L. Godo, P. Hájek and  M. Navara proved that $Alg^\ast {\cal G}_{\s}$ is the class of Gödel algebras that satisfy the equations of the form $\varphi \approx 1$ where $\varphi$ ranges over the  axiom schemata $(\s 1)$, $(\s 2)$, $(\s 3)$, $(\Delta 1)$, $(\Delta 2)$ and $(\Delta 5)$. Note that, by Lemma \ref{prop aritmeticas de $G_sim$} and \cite[Lemma 4]{EG},  these algebras are precisely the $G_{\s}$-algebras defined in this article.

 On the other hand, D. Castaño et al. proved that $Alg^\ast {\cal S}5'(BL)$ is the class of  monadic $BL$-algebras.  Taking these results into account,  and the fact that the axioms $(\s 2)$, $(\lozenge 1)$ and $(\lozenge 2)$  correspond to the conditions $(N)$, $(Q)$ and $(C)$, the proof of this result is immediate.
\end{Proof}

Now, using  (\ref{semantica algebraica equivalente}) for the case ${\cal L}= {\cal S}5({\cal G}_{\s})$  we have the following completeness  result.

\begin{Theorem} 
For every formula $\varphi$ and set of formulas $\Gamma$ in ${\cal S}5({\cal G}_{\s})$:
$$\Gamma \vdash_{{\cal S}5({\cal G}_{\s})} \varphi \mbox{ \ \  if and only if \ \ } \Gamma \models_{\mathbb{CMG}_{\s}} \varphi,$$
\end{Theorem}

We next prove a result connecting completeness with functional algebras similar to the one given in \cite{Complet}.

\begin{Theorem} \label{completitud via funcionales}
The following conditions are equivalent:

\begin{enumerate}[$(1)$]
\item  For any formula $\varphi$ and any set of formulas  $\Gamma$, we have $$\Gamma \vdash_{{\cal S}5({\cal G}_{\s})} \varphi \mbox{ \ \  if and only if \ \ } \Gamma \models_{{\cal S}5({\mathbb G}_{\s})} \varphi$$

\item The variety $\mathbb{CMG}_{\s}$ is generated, as a quasivariety, by its functional monadic $G_{\s}$-algebras.

\end{enumerate}

\end{Theorem}
\begin{Proof}
The proof runs in a similar way to the corresponding proof in \cite[Theorem 2.1]{Complet} taking into account that  ${\cal S}5({\mathbb G}_{\s})$ is  finitary because it is equivalent to the one-variable monadic fragment  of the predicate calculus ${\cal G}\forall_{\s}$ and this latter is finitary, and $\mathbb{CMG}_{\s}$ is the equivalent algebraic semantics of the logic ${\cal S}5({\cal G}_{\s})$ (Proposition \ref{completitud fuerte con S5 y CMG}).
\end{Proof}

Finally, since $(2)$ in  Theorem \ref{completitud via funcionales} holds (see Corollary \ref{CMG generada como cuaisvariedad por funcionales con punto fijo}) we have the following completeness result    for ${\cal S}5({\cal G}_{\s})$.
 
\begin{Theorem}[Completeness]

For any formula $\varphi$ and any set of formulas  $\Gamma$, it follows that

$$\Gamma \vdash_{{\cal S}5({\cal G}_{\s})} \varphi \mbox{ \ \ if and only if \ \ } \Gamma \models_{{\cal S}5({\mathbb G}_{\s})} \varphi \mbox{ \ \  if and only if \ \ } \Gamma \models_{\mathbb{CMG}_{\s}}\varphi. $$

\end{Theorem}

Thus, we have characterized the equivalent algebraic semantics for the one-variable monadic fragment of the first-order logic ${\cal G}\forall_\sim$ as the variety $\mathbb{CMG}_\sim$.  As a consequence, there exists a correspondence between the axiomatic extensions of this calculus and the subvarieties of $\mathbb{CMG}_\sim$. For this reason, it would be important to describe all subvarieties of this variety. However, considering the complexity of this problem, we believe that a way to proceed with our work would be to intend to resolve this problem for the subvariety of $\mathbb{CMG}_\sim$ generated by totally ordered algebras. 

To carry this out, we believe it would be very useful to use the results obtained in \cite{Laura} from the variety $\mathbb{MG}$ and to obtain a suitable topological duality from $\mathbb{CMG}_\sim$ that allows us to deepen the study of this variety.


\bigskip

\bigskip

{\sc Diego Castaño}, {\tt diego.castano@uns.edu.ar}

Departamento de Matemática, Universidad Nacional del Sur (UNS), Bahía Blanca, Argentina

Instituto de Matemática (INMABB), Universidad Nacional del Sur (UNS)-CONICET, Bahía Blanca, Argentina

\bigskip

{\sc Valeria Castaño}, {\tt cvaleria@gmail.com}

Departamento de Matemática, Universidad Nacional del Comahue, Neuquén, Argentina

IITCI, UNCo-CONICET, Neuquén, Argentina

\bigskip

{\sc José Patricio Díaz Varela}, {\tt usdiavar@criba.edu.ar}

Departamento de Matemática, Universidad Nacional del Sur (UNS), Bahía Blanca, Argentina

Instituto de Matemática (INMABB), Universidad Nacional del Sur (UNS)-CONICET, Bahía Blanca, Argentina

\bigskip

{\sc Marcela Muñoz Santis}, {\tt santis.marcela@gmail.com}

Departamento de Matemática, Universidad Nacional del Comahue, Neuquén, Argentina

\end{document}